\documentclass[preprint,review,12pt]{elsarticle}

\makeatletter
\def\ps@pprintTitle{%
  \let\@oddhead\@empty
  \let\@evenhead\@empty
  \def\@oddfoot{\reset@font\hfil\thepage\hfil}
  \let\@evenfoot\@oddfoot
}
\makeatother

\usepackage[utf8]{inputenc}
\usepackage{fullpage,times}
\usepackage{amsmath,amssymb}
\usepackage{graphicx}
\usepackage{algorithm}
\usepackage{algorithmicx}
\usepackage{subfig}
\usepackage{caption}
\usepackage{url}
\usepackage[numbers]{natbib} 
\usepackage{placeins}
\usepackage{siunitx}
\usepackage{url}
\usepackage{xspace}
\usepackage[table,xcdraw]{xcolor}

\usepackage{hyperref}
\hypersetup{
colorlinks={true},
linkcolor={Blue3},
urlcolor={Blue3},
citecolor={Blue3},
pdfstartview={FitH},
pdfauthor={Author},
pdftitle={Title},
pdfsubject={Subject}
plainpages=false
}

\usepackage{enumitem}

\usepackage{booktabs} 
\usepackage{multirow}
\usepackage{multicol}
\usepackage{makecell}

\usepackage[toc,page]{appendix}

\usepackage{setspace}

\definecolor{wisconsin-red}{rgb}{0.6,0,0}
\definecolor{darkgreen}{rgb}{0.2,0.6,0.2}
\definecolor{maroon}{rgb}{0.5, 0.0, 0.0}
\definecolor{violet}{rgb}{0.75, 0.0, 1.0}
\definecolor{lightgray}{gray}{0.9}
\definecolor{navyblue}{rgb}{0.0, 0.0, 0.5}
\definecolor{darkmidnightblue}{rgb}{0.0, 0.2, 0.4}
\definecolor{midnightblue}{rgb}{0.0,0.4,0.85}
\definecolor{Gray}{gray}{0.75}
\definecolor{darkgreen}{rgb}{0,0.5,0}
\definecolor{apricot}{rgb}{0.98, 0.81, 0.69}

\newcolumntype{C}[1]{>{\centering\arraybackslash}p{#1}}
\newcolumntype{P}[1]{>{\raggedright\arraybackslash}p{#1}}
\newcolumntype{L}[1]{>{\raggedleft\arraybackslash}p{#1}}

\usepackage{upgreek}

\newcommand{\NeurADP}{\texttt{NeurADP}}
\newcommand{\NeurADPFixed}{\texttt{NeurADP-Fixed}}
\newcommand{\Myopic}{\texttt{Myopic}}
\newcommand{\Direct}{\texttt{Direct}}
\newcommand{\DRL}{\texttt{DRL}}
\newcommand{\Time}{\mathcal{T}}
\newcommand{\GraphNetwork}{\mathcal{N}}

\newcommand{\BatchingSet}{B}
\newcommand{\FeasibilitySet}{F}

\newcommand{\DecisionSet}{\mathbf{A}}

\newcommand{\CouriersState}{C}
\newcommand{\OrdersState}{O}
\newcommand{\CouriersSet}{\mathcal{C}}
\newcommand{\OrdersSet}{\mathcal{O}}
\newcommand{\TimeIntervals}{\texttt{$\delta$}}
\newcommand{\Locations}{\texttt{$\mathcal{L}$}}
\newcommand{\Edges}{\texttt{$\mathcal{E}$}}
\newcommand{\CurrentTime}{\texttt{$t$}}
\newcommand{\StateNote}{\texttt{$S$}}

\newcommand{\RequestDeadline}{\texttt{$dead$}}
\newcommand{\TravelTime}{\texttt{time}}

\newcommand{\MaxDelay}{\texttt{delay\textsubscript{\texttt{max}}}}

\newcommand{\MovementPoint}{f}

\newcommand{\ExogenousInformation}{W}
\newcommand{\NullAction}{\varnothing}
\newcommand{\DecisionTuple}{\textbf{a}}
\newcommand{\PostDecisionState}{S^{\texttt{Courier-Post}}}
\newcommand{\FirstTransition}{{\texttt{statepost}}}
\newcommand{\SecondTransition}{{\texttt{statenext}}}

\newcommand{\IndividualVehicleDropOff}{\texttt{$drop$}}
\newcommand{\Reward}{R}
\newcommand{\ImmediateReward}{r}
\newcommand{\PostDecision}{\texttt{Post}}

\newcommand{\EstimatedV}{\bar{v}}

\newcommand{\StepSize}{\alpha}

\newcommand{\SingleDay}{n}
\newcommand{\MarginalValue}{\hat{v}}

\newcommand{\QueueCadrinality}{D}

\newcommand{\IndividualCourier}{\texttt{$c$}}
\newcommand{\IndividualOrder}{\texttt{$o$}}
\newcommand{\CourierShiftStart}{\texttt{shift}}
\newcommand{\CourierReturn}{\texttt{ret}}
\newcommand{\CourierOrders}{\texttt{ords}}
\newcommand{\OrderDestination}{\texttt{dest}}
\newcommand{\OrderDeadline}{\texttt{dead}}

\newcommand{\CourierPost}{\texttt{Courier-Post}}
\newcommand{\FeasibleBatching}{b}
\newcommand{\MaxQueueSize}{\texttt{queue\textsubscript{\texttt{max}}}}
\newcommand{\ShiftLength}{\texttt{shift$_{\texttt{length}}$}}

\newcommand{\IndividualReward}{\texttt{$q$}}
\newcommand{\PParameter}{\texttt{$\beta$}}
\newcommand{\ActionDelay}{\texttt{$\omega$}}

\newcommand{\PermutationSet}{\mathcal{Z}}
\newcommand{\SinglePermutation}{\sigma}
\newcommand{\SingleIndex}{d}

\newcommand{\MatchingVariable}{a}

\begin{document}\sloppy
\setlength{\parindent}{2em}

\begin{frontmatter}
\title{Neural Approximate Dynamic Programming for the Ultra-fast Order Dispatching Problem} 

\author[1]{Arash Dehghan}
\ead{arash.dehghan@torontomu.ca}

\author[1]{Mucahit Cevik\corref{cor1}%
\fnref{fn1}}
\ead{mcevik@torontomu.ca}

\author[2]{Merve Bodur}
\ead{merve.bodur@ed.ac.uk}

\cortext[cor1]{Corresponding author}
\fntext[fn1]{Toronto Metropolitan University, Toronto, ON, Canada}
\address[1]{Toronto Metropolitan University, Toronto, ON, Canada}
\address[2]{University of Edinburgh, Edinburgh, UK}

\begin{abstract}
Same-Day Delivery (SDD) services aim to maximize the fulfillment of online orders while minimizing delivery delays but are beset by operational uncertainties such as those in order volumes and courier planning. 
Our work aims to enhance the operational efficiency of SDD by focusing on the ultra-fast Order Dispatching Problem (ODP), which involves matching and dispatching orders to couriers within a centralized warehouse setting, and completing the delivery within a strict timeline (e.g., within minutes). 
We introduce important extensions to ultra-fast ODP such as order batching and explicit courier assignments to provide a more realistic representation of dispatching operations and improve delivery efficiency. As a solution method, we primarily focus on NeurADP, a methodology that combines Approximate Dynamic Programming (ADP) and Deep Reinforcement Learning (DRL), and our work constitutes the first application of NeurADP outside of the ride-pool matching problem. NeurADP is particularly suitable for ultra-fast ODP as it addresses complex one-to-many matching and routing intricacies through a neural network-based VFA that captures high-dimensional problem dynamics without requiring manual feature engineering as in generic ADP methods.
We test our proposed approach using four distinct realistic datasets tailored for ODP and compare the performance of NeurADP against myopic and DRL baselines by also making use of non-trivial bounds to assess the quality of the policies.
Our numerical results indicate that the inclusion of order batching and courier queues enhances the efficiency of delivery operations and that NeurADP significantly outperforms other methods. Detailed sensitivity analysis with important parameters confirms the robustness of NeurADP under different scenarios, including variations in courier numbers, spatial setup, vehicle capacity, and permitted delay time. 
\end{abstract}
\end{frontmatter}
\vspace{-8pt}
\section{Introduction}
The widespread adoption of online shopping, particularly accelerated by the COVID-19 pandemic, has transformed traditional markets in recent years and compelled many businesses to embrace streamlined direct delivery of products to customers \citep{dayarian2020same}. One notable consequence of this shift is the emergence of Same-Day Delivery (SDD) services, which have fundamentally changed shopping behaviors by offering the convenience of online ordering and near-instant access to products. The SDD has seen a remarkable growth in recent years, with a valuation of \$5.77 billion in the United States in 2019 and a projected value of \$15.6 billion by 2024 \citep{statista}. Recognizing the evolving dynamics of the retail landscape, major players such as Target, Walmart, and Amazon have all acknowledged the significance of providing competitive same-day shipping options and have looked to expand their same-day shipping services \citep{thomas2019target}. As a result of this rapid expansion, centralized warehouses have become the central hub for managing incoming online orders and dispatching fleets of couriers, all with the goal of providing efficient and prompt service.

With the growing popularity of SDD operations, it is crucial to prioritize operational efficiency. The primary goal of SDD operations is to maximize the fulfillment of online orders while minimizing delivery delays. Nevertheless, SDD services naturally encompass several considerations that need to be taken into account in managing delivery operations. The courier shift schedules, vehicle capacities and dynamic routing of delivery couriers are some of the important considerations in this regard. Additionally, SDD operations involve a substantial level of uncertainty that stems from factors such as the timing, volume, deadlines, and destination locations of orders. These multitude of factors pose significant challenges for SDD operators, who must navigate them to provide efficient services. To this effect, \citet{voccia2019same} introduce the Same-Day Delivery Problem (SDDP) as a framework to define the complex decision-making and routing logistics involved in ensuring the timely delivery of online orders within strict time constraints. 

The SDDP may be decomposed into two distinct sub-problems: the Vehicle Routing Problem (VRP) and the Order Dispatching Problem (ODP). The VRP addresses the routing aspect, while the ODP concentrates on the matching and dispatching components of the problem. Both sub-problems are usually relevant in the context of a centralized warehouse handling stochastic order arrivals for cost-effective and timely delivery to customers. The VRP involves minimizing total vehicle travel distance, time, or cost while accounting for congestion, capacity, and time windows. On the other hand, the ODP involves assigning orders to couriers, minimizing fulfillment time and avoiding capacity breaches, factoring in location, time windows, and order size. 

In this paper, we focus on the ultra-fast ODP, which involves the matching and dispatching aspects of the SDDP. Specifically, we explore a centralized decision-making problem in which a warehouse dispatches fleets of couriers, each with their own shift schedules, to maximize the number of orders served throughout the day. These orders arrive stochastically, and the warehouse's primary objective is to ensure \textit{ultra-fast} deliveries, e.g., completing them within minutes. This requirement of urgent delivery introduces a critical time constraint, distinguishing it from other SDDP/ODP works which allow for more lenient delivery timelines. Furthermore, despite the widespread adoption of this rapid delivery approach by global delivery giants such as \textit{Getir}, a renowned Turkish delivery service, and \textit{Gorillas}, a Germany-based platform specializing in swift grocery and essential item deliveries, this particular delivery setting is not well-studied in the literature. In this regard, our paper contributes to the existing literature on the SDDP and ODP by building upon the work of \citet{kavuk2022order}, which focuses on an order dispatching problem based on Getir's operations. 
Specifically, they develop a deep reinforcement learning (DRL) approach for ODP to make informed decisions regarding only the acceptance or rejection of incoming orders to a single depot while the order assignments to the couriers are based on predefined rules. The main contributions of our paper are summarized in what follows.

We propose a novel Markov decision process (MDP) model that introduces several innovative features and capabilities to the single-depot ODP. In particular, in terms of operational enhancements (1) we employ \textit{batching} as a means of enhancing the efficiency of the order dispatching operations and we utilize \textit{courier queues} which enable the concurrent handling and emptying of all orders within the queue, rather than being limited to serving one order at a time, (2) we make \textit{explicit courier assignments} to optimize the allocation of orders to couriers, and (3) we enforce \textit{hard deadlines} to ensure timely delivery of orders.
By incorporating these considerations, our proposed MDP model for ODP provides a more comprehensive and realistic representation of the dispatching process, thereby enhancing its practical relevance.

Inspired by its effective application to the ride-pool matching problem, and observing the suitability of our problem's structure to leverage its strengths, we adopt the Neural Approximate Dynamic Programming (NeurADP) as the solution approach. NeurADP is an innovative methodology introduced by \citet{shah2020neural} which combines Approximate Dynamic Programming (ADP) and DRL techniques and has exclusively been applied within the ride-sharing framework. Our study constitutes the first application of NeurADP beyond its original context, expanding its potential applications and further demonstrating its effectiveness in addressing real-world dynamic decision-making problems. We note that NeurADP is well-suited for ultra-fast ODP as it can address complex one-to-many matching and routing intricacies through a neural network-based value function approximation (VFA) that captures high-dimensional problem dynamics without requiring manual feature construction as in many other ADP frameworks.
In order to demonstrate the effectiveness of NeurADP, we compare it with a large set of myopic and DRL baselines. We also conduct a sensitivity analysis to investigate the influence of various factors on the performance of NeurADP. 

In our numerical study, to support our research and facilitate comprehensive evaluations, we introduce three original datasets specifically tailored for order dispatching operations in addition to considering a commonly used dataset from the literature. These datasets capture diverse real-world scenarios and provide a rich environment for training and testing our proposed methods. The availability of these datasets benefits future researchers in the field, fostering further advancements in the study of ODP and SDDP.
Furthermore, we explore the effects of the number of agents considered, the spatial setup, the allowed vehicle capacity, and the permitted delay time. This analysis enhances our understanding of the robustness and adaptability of NeurADP in diverse scenarios, offering valuable insights to practitioners seeking to implement this approach in real-world applications.
Importantly, our analysis also provides several managerial insights related to the ultra-fast ODP:
\begin{itemize}
    
    \item NeurADP improves order fulfillment by 6.7\%- 16.9\% in the baseline configuration compared to the benchmark policies thanks to its ability to intelligently batch orders and optimize courier utilization. Furthermore, the greatest benefits of NeurADP-based policies are observed when there are fewer couriers who are working at or near full capacity or when the operation faces tighter delivery schedules.

    \item In terms of fulfillment strategy, NeurADP policy shows that the companies can enhance their efficiency by strategically batching the orders and ensuring the swift return of couriers to the warehouse, rather than attempting to maximize the number of orders loaded onto each courier's vehicle. That is, the companies can potentially serve more orders with quicker turnaround which could be more beneficial than simply loading couriers with as many orders as possible, especially in a high-demand and fast-paced delivery environment.

    \item The significance of policy intelligence in optimizing order fulfillment varies depending on the operation environment, particularly based on the sparsity of the delivery locations and their distance to the central warehouse. In cases where deliveries are more dispersed, the effectiveness of a policy becomes markedly more crucial. On the other hand, when delivery locations are densely clustered and closer to the warehouse, the relative importance of having an intelligent policy diminishes. 
    This suggests that companies operating in diverse geographic settings should tailor their dispatching policies to the specific distribution characteristics of each area to optimize courier efficiency and order fulfillment rates.
\end{itemize}

The remainder of the paper is organized as follows. Section~\ref{literaturereview} offers a comprehensive review of the relevant SDD literature in the context of both the VRP and ODP, better positioning our research within the existing body of work. Section~\ref{problemformulation} provides a formal description of the problem setting for our ODP. In Section~\ref{solutionmethodology}, we describe the NeurADP solution methodology. Details regarding the datasets and benchmark policies used in the experiments are provided in Section~\ref{experimentalsetup}. The results of the computational experiments are presented in Section~\ref{results}, followed by a conclusion in Section~\ref{conclusion} that summarizes the research findings and suggests avenues for future research.

\section{Literature Review} \label{literaturereview}
We review the relevant literature by exploring the challenges and solutions related to the SDDP and its sub-problems, particularly ODP, highlighting the latest research and findings in this field, as well as placing our work and contributions within the broader context of the existing literature. In order to tackle the challenges in the SDDP, various strategies have been proposed, encompassing heuristic algorithms, machine learning models, combinatorial optimization models, reinforcement learning (RL) algorithms, and market-based mechanisms. Table~\ref{table:references} presents the most relevant studies to ours from the SDDP literature. This table comprises six indicators regarding the problem context and solution methodology. These are ``Solution Technique'', which describes the approach employed to solve the problem, ``Large Capacity'', which denotes whether couriers are allowed to carry multiple sets of orders simultaneously, ``Multi-Courier'', which indicates consideration of more than one courier, ``Hard Deadlines'', which signifies the presence of strict delivery deadlines, ``Bundling'', which pertains to the possibility of bundling/batching orders at the same decision-making step, and ``Shifts'', which indicates the incorporation of courier shifts into the respective problem formulation.


\setlength{\tabcolsep}{2.5pt} 
\renewcommand{\arraystretch}{1.3} 
\begin{table}[!ht]
\centering
\caption{Summary of relevant studies. (ADP: Approximate Dynamic Programming, DRL: Deep Reinforcement Learning, CH: Combinatorial Heuristic, MIP: Mixed Integer Programming, PFA: Policy Function Approximation, NeurADP: Neural Approximate Dynamic Programming)}\label{table:references}
\resizebox{0.99\textwidth}{!}{
\begin{tabular}{P{0.29\textwidth}C{0.13\textwidth}C{0.13\textwidth}C{0.13\textwidth}C{0.13\textwidth}C{0.13\textwidth}C{0.13\textwidth}C{0.13\textwidth}} 
\toprule
\textbf{Study} & \textbf{Problem} & \textbf{Solution Technique} & \textbf{Large Capacity} & \textbf{Multi-Courier} & \textbf{Hard Deadlines} & \textbf{Bundling} & \textbf{Shifts} \\ 
\midrule
\citet{ulmer2019ADP} & VRP  & ADP & \checkmark & & & \checkmark & \\
\citet{joe2020deep} & VRP & DRL &  \checkmark & \checkmark & & & \\
\citet{cote2021dynamic} & VRP & CH & \checkmark & \checkmark &    \checkmark & & \\
\citet{ngu2022decentralised} & VRP & DRL &  \checkmark & \checkmark & & & \\
\citet{klapp2017} & ODP & MIP & \checkmark & & &           \checkmark & \\
\citet{kavuk2022order} & ODP & DRL &                            & \checkmark & & & \checkmark \\
\citet{ulmerstreng2019} & ODP & PFA & \checkmark & \checkmark & \checkmark & \checkmark & \\
\citet{ulmerdrone2019ADP} & ODP & PFA & \checkmark & \checkmark & \checkmark & & \\
\citet{chen2022deep} & ODP & DRL & \checkmark & \checkmark & \checkmark & & \\
\citet{cardona2022same} & ODP & PFA & \checkmark & \checkmark & \checkmark & & \\ 
\midrule
\textbf{Our Work} & ODP & NeurADP & \checkmark & \checkmark & \checkmark & \checkmark & \checkmark \\ 
\bottomrule
\end{tabular}
}
\end{table}

VRP has been explored within the context of the SSDP in various works. \citet{ulmer2019ADP} presented an ADP-based order assignment/dispatching and routing policy which allows same-day delivery vehicles to better integrate dynamic requests into delivery routes through preemptive depot returns. \citet{joe2020deep} combined DRL with a simulated annealing-based routing heuristic for a dynamic VRP. Their method uses a state representation based on the total cost of the remaining routes of the vehicles. Additionally, \citet{cote2021dynamic} proposed a re-optimization heuristic and a branch-and-regret heuristic that uses sampled scenarios to anticipate future events to address a variation of the VRP that involves urgent deliveries of time-sensitive orders. \citet{ngu2022decentralised} presented a decentralized multi-agent RL approach in formulating and solving the VRP using a parameter-sharing deep Q-network. \citet{ulmer2020dynamic} proposed a method called anticipatory pricing and routing policy to improve the cost-efficiency of same-day delivery for e-commerce retailers. By dynamically adjusting prices based on delivery deadlines and using a guided offline VFA, this policy incentivizes customers to select efficient delivery options, allowing the fleet to serve more orders and increase revenue. Finally, \citet{dayarian2020same} explored the concept of drone replenishment in the context of same-day home delivery. They introduced the VRP with drone resupply and proposed various algorithms to optimize the delivery process, quantifying the potential benefits of using drones for delivery vehicles.

Several studies have specifically focused on the dispatching aspect of the SDDP. \citet{cardona2022same} explored a two-echelon fleet approach that utilizes intra-route replenishment and policy function approximation (PFA) based on real-life geographical distributions to optimize fleet configuration and maintain service levels. Whereas \citet{ulmerdrone2019ADP} used fleets of heterogeneous drones and vehicles to perform deliveries, utilizing PFA based on geographical districting to decide which delivery method of transportation to use. Similarly, \citet{chen2022deep} proposed a same-day delivery system using both vehicles and drones and presented a deep Q-learning approach to learn the value of assigning customer orders to either drones, vehicles, or not offering service at all. To enable real-time dispatch decisions that balance speedy delivery with consolidation, \citet{ulmerstreng2019} introduced a novel same-day delivery approach that combines autonomous vehicles with pickup stations and utilized a PFA approach. In addition, \citet{klapp2017} formulated an arc-based Mixed Integer Programming (MIP) model and designed local search heuristics to solve the deterministic version of the ODP and derived an apriori solution for the stochastic case. Finally, \citet{kavuk2022order} presented a DRL approach to solving the ODP for ultra-fast delivery, using deep Q-networks to learn the actions of warehouses and considering two reward functions: one related to the number of orders served and the other to minimize delivery delays.


In this paper, we consider the matching and dispatching components of the SDDP, particularly ODP, differentiating our focus from prior research that primarily concentrated on the routing elements \citep{cote2021dynamic, joe2020deep, ngu2022decentralised, ulmer2019ADP}. In particular, we consider a centralized decision-making scenario where a warehouse is tasked with coordinating fleets of couriers, each operating on its own shift schedule in order to optimize the total number of orders fulfilled throughout the day. These orders are subject to stochastic arrivals, and the primary goal of the warehouse is to achieve ultra-fast deliveries, aiming to complete them within a matter of minutes. This imperative for rapid delivery imposes a crucial time constraint, setting it apart from previous SDDP and ODP works which permit more relaxed delivery time-frames \citep{joe2020deep, klapp2017, ngu2022decentralised, ulmer2019ADP}. The imperative to handle deliveries within minutes necessitates the capability of real-time decision-making, thereby challenging the feasibility of past traditional offline solutions \citep{klapp2017}. The dynamic nature of the problem is further accentuated by the stochastic arrival of orders, requiring dispatching and matching strategies to flexibly adapt to varying patterns, distinguishing it from the approaches seen in previous works \citep{restrepo2019integrated, ulmer2019ADP}. Moreover, the added complexity arises from the individual courier shift schedules, which demand meticulous coordination and optimization to meet the stringent ultra-fast delivery criteria. This inclusion of courier shifts sets our work apart from prior studies such as \citep{cardona2022same, chen2022deep, ulmerstreng2019, ulmerdrone2019ADP} which did not encompass this facet. In the pursuit of achieving such rapid deliveries, the allocation of couriers to orders becomes a task of precise resource allocation optimization, a contrast to the focus of traditional ODP studies, which primarily revolve around dispatching timing and transportation mode decisions ( e.g., see~\citep{ cardona2022same, ulmerstreng2019, ulmerdrone2019ADP}). Despite the widespread adoption of this rapid delivery approach by delivery corporations such as \textit{Getir} and \textit{Gorillas}, this particular delivery setting has not been well-studied in the existing literature.

We introduce a comprehensive set of improvements for order dispatching operations in the context of ultra-fast delivery, building upon the work of \citet{kavuk2022order}. Their research closely aligns with ours, particularly as they focus on addressing the ODP encountered by \textit{Getir}. In their study, \citet{kavuk2022order} employed DRL to determine only the acceptance or rejection decisions for incoming orders. However, their work does not take into account important ODP considerations including courier assignment and order batching, and importantly, their framework does not impose strict delivery deadlines, rather penalizing the delays. Furthermore, their empirical analysis is limited to a single dataset and the rule-based baselines for comparative analysis. To further enhance the existing problem framework, our work introduces innovative features and capabilities, namely, order batching, explicit courier assignment and hard deadlines. Hence, it helps streamline the coordination and efficiency of order dispatching, contributing to the practical relevance of this problem. Moreover, we extend the scope of the problem to encompass different urban settings and larger-scale dispatching operations, involving more agents, orders, and a broader geographical area. This expansion enables us to capture the intricacies and challenges of managing substantial dispatching tasks, providing valuable insights for real-world applications. To tackle these challenges, we adapt the innovative NeurADP approach for order dispatching, which was originally designed for ride-sharing, marking its first application outside its original context and showcasing its effectiveness in dynamic decision-making problems. Moreover, to support our research and facilitate comprehensive evaluations, we introduce three novel tailored datasets for order dispatching. 

\section{Problem Description and Formulation} \label{problemformulation}
We present a dynamic order dispatching model that aims to efficiently match couriers with incoming batches of online orders. Our model considers the spatial and temporal demand patterns of the orders, which arrive dynamically over a 24-hour decision horizon and are served by a single centralized depot. This choice of a single central depot is particularly important in the context of ultra-fast delivery, where efficiency and speed are paramount. For instance, the logistics of coordinating multiple depots may introduce unnecessary delays and complexities, ultimately hindering the goal of rapid order fulfillment. Orders are generated stochastically and have specific delivery deadlines based on their arrival time. Once an order is assigned to a courier, it is accepted into the system and its delivery prior to its designated drop-off deadline is guaranteed. Moreover, the model takes into account a predetermined group of heterogeneous couriers available during the planning horizon, considering their capacity constraints and shift schedules. All couriers have individual shift start times, with each shift lasting six hours (excluding breaks) to reflect real-world shift lengths.

Our model incorporates several key problem specifications in the ultra-fast delivery setting to efficiently manage the dispatch and delivery process. First, multiple orders may arrive at any decision epoch, and couriers located at the warehouse are promptly dispatched upon being matched to these orders. Secondly, once a courier is dispatched with a set of orders, they must complete all assigned deliveries before returning to the depot, precluding any preemptive returns. To aid in this process, a queue is maintained for each courier with a capacity equivalent to their vehicle's limit. This queue accommodates both pending orders awaiting pickup and delivery, as well as new orders which may be matched to an on-shift courier as they continue their deliveries. Orders are incorporated into the queue only if their inclusion maintains adherence to constraints regarding timely delivery of all orders within the queue, steering clear of overloading courier vehicles, and staying within courier shift duration during order deliveries. Furthermore, the queue of orders is rearranged prior to the courier dispatching from the depot so as to optimize the route from the warehouse to each order destination location and back. However, once an order is assigned to a specific courier's queue, it cannot be transferred to another courier's queue. Lastly, in line with prior research on ODP, unmatched orders beyond their arrival period are assumed to exit the system~\citep{kavuk2022order}. This assumption reflects customers' general expectation of timely confirmation regarding the acceptance of their requests. 

The primary objective of our model is to maximize the total number of online orders fulfilled within the decision horizon. To achieve this, our assignment decisions consider future order arrival uncertainties and the potential downstream impact of current decisions. To handle the complexity of these decisions, we formulate an MDP model and adapt a NeurADP solution framework, enabling effective real-time decision-making under uncertainty. To this end, we partition the finite planning horizon into discrete time intervals, with each interval having a duration of \TimeIntervals~(e.g., five minutes). We assume that decisions are made at the onset of each interval, while exogenous information is observed continuously throughout. Following each interval, the state of the system is updated by incorporating the decisions and the observed external information. The collection of epochs for decision-making is denoted as $\Time:=\{0,\ldots,T\}$. At each decision epoch, the aim is to match ``available'' couriers with incoming orders. The availability of a courier is determined by several factors, including whether they are on their shift, their available capacity, and whether adding a new order to their assigned order set would comply with the maximum allowed delay for any order and would not extend the courier beyond their shift end time. Both the couriers who are stationed at the warehouse, as well as those who are away from the warehouse making deliveries, are eligible to be paired with incoming orders, provided they satisfy the availability constraints. Once paired with a batch of incoming orders, couriers located at the warehouse are promptly dispatched so as to adhere to the ultra-fast delivery requirements, while those making deliveries maintain a queue for orders to pick up and deliver after completing their ongoing assignments. The courier's queue accepts orders up until the moment the courier returns to the warehouse, at which point it is emptied, and the courier is promptly dispatched with any orders that may have accumulated in their queue. Couriers who are off-shift cannot be matched with orders, keeping their queues empty.

Each courier has a predetermined start time for their shift, with each shift spanning a duration of six hours, without any scheduled breaks. 
Incoming orders are typically associated with a specific delivery deadline which can be set in different ways depending on the company policies.
For instance, in their paper, \citet{kavuk2022order} consider a 45-minute delivery time for any given order and employ a reward function that promotes fast deliveries below this 45-minute target.
Lastly, it should be noted that, at each time step, orders have the potential to be consolidated (i.e., batched) either with other concurrent orders or with previously assigned orders for each respective courier. The network housing the depot is characterized as $\GraphNetwork=(\Locations,\Edges)$, where $\Locations=\{0,1,\ldots,L\}$ corresponds to the depot and customer locations, and $\Edges=\Locations \times \Locations$ represents the distance between each respective pair of locations, determined by the Haversine distances. Given two locations $\ell$, $\ell'$ $\in \Locations$, we denote the travel time between $\ell$ and $\ell'$, leaving $\ell$ at time $\CurrentTime$, by $\TravelTime_{\CurrentTime}(\ell,\ell')$, such that $\ell=0$ always denotes the depot location. Next, we present the components that make up the MDP model.

\subsection{State Variables}
The state of the system at time $\CurrentTime \in \Time$ is defined by $\StateNote_\CurrentTime=(\CouriersState_\CurrentTime,\OrdersState_\CurrentTime)$, such that $\CouriersState_\CurrentTime$ represents the state of all couriers, and $\OrdersState_\CurrentTime$ the state of all incoming orders awaiting delivery. The state of an individual courier may be represented as a three-dimensional attribute vector defined by $\IndividualCourier=(\IndividualCourier_\CourierShiftStart,\IndividualCourier_\CourierReturn, \IndividualCourier_\CourierOrders)\in \CouriersSet$ with $\CouriersSet$ denoting the set of possible courier states. In this representation, $\IndividualCourier_\CourierShiftStart$ indicates the time at which the shift of the courier starts. Moreover, if the courier is away from the warehouse fulfilling deliveries, $\IndividualCourier_\CourierReturn$ signifies the time required for them to complete their deliveries and return to the warehouse. This value is set to zero if the courier is not on shift or is already at the warehouse. Lastly, $\IndividualCourier_\CourierOrders$ represents the courier's queue and encompasses significant details about the courier's current tasks, including the orders they are currently assigned for delivery upon their return to the warehouse. The sequence of online orders assigned to the queue is optimized to minimize travel time for delivering all orders, and it is rearranged prior to a courier departing from the depot to make deliveries. Subsequently, the state of an online order is represented by a two-dimensional attribute vector denoted as $\IndividualOrder=(\IndividualOrder_\OrderDestination, \IndividualOrder_\OrderDeadline) \in \OrdersSet$ with $\OrdersSet$ denoting the set of possible incoming orders. Here, $\IndividualOrder_\OrderDestination$ denotes the destination of the order, while $\IndividualOrder_\OrderDeadline$ corresponds to the specific delivery deadline time. When an order is received between decision epochs $t-1$ and $t$, the deadline attribute $\IndividualOrder_\OrderDeadline$ is determined at the start of epoch $t$ using the equation
\begin{equation}
    \IndividualOrder_\OrderDeadline = \CurrentTime + \TravelTime_\CurrentTime(0, \IndividualOrder_\OrderDestination) + \MaxDelay. 
    \label{eq:OrdDeadline}
\end{equation}
where \MaxDelay indicates the maximum allowed time beyond the original travel duration from the depot to the order's drop-off location. While this method of calculating the delivery deadlines is slightly different than \citet{kavuk2022order}'s approach, we note that our proposed framework can accommodate alternative ways of setting the delivery deadlines. Whereas we note that explicitly setting a delay parameter can help setting more realistic customer expectations.

\subsection{Decision Variables}\label{sec:DecisionVariables}
At each decision epoch $\CurrentTime \in \Time$, we determine the matching between available couriers (i.e., those idly waiting in the warehouse or the busy couriers with available space in their queues) and incoming online orders considering the current system state. More specifically, we begin by examining the feasibility of grouping the set of incoming orders into batches, taking into account the order drop-off deadlines, and then evaluate the potential for assigning a specific batch to a courier, taking into consideration both capacity and timing limitations. To evaluate the feasibility of batching a set of orders together, or taking a single order by itself, we make sure whether a batch is able to be delivered before each order's respective drop-off deadline. In other words, a batching is feasible if there exists a viable route for a courier to deliver each order in the batch, as well as the orders it is currently assigned to, prior to each order's respective deadline beginning from the warehouse. Furthermore, when deciding if a courier can be paired with a batch of orders, we consider the following factors: (i) the courier's active status and current shift, (ii) whether the newly assigned batch pushes the courier's queue beyond the allowed limit, (iii) the courier's ability to deliver all orders before their specific deadlines, and (iv) whether the courier can complete all deliveries and return to the depot before their shift ends. We subsequently define the collection of actions taken at time $\CurrentTime$ as $\DecisionTuple_\CurrentTime \in \DecisionSet_\CurrentTime(\StateNote_\CurrentTime)$, such that $\DecisionSet_\CurrentTime(\StateNote_\CurrentTime)$ denotes the set of all feasible actions for state $\StateNote_\CurrentTime$. These actions encompass both the matching of couriers to order batches, as well as the determination of delivery sequencing within each courier's assigned orders. 

We define the reward collected at time step $\CurrentTime \in \Time$ as follows:
\begin{equation}
    \Reward_\CurrentTime (\DecisionTuple_\CurrentTime)=
    \sum_{\IndividualCourier\in \CouriersState_\CurrentTime} \left(\PParameter \cdot \IndividualReward_\CurrentTime(\DecisionTuple_{\CurrentTime \IndividualCourier})~-~\ActionDelay_{\CurrentTime}(\DecisionTuple_{\CurrentTime \IndividualCourier})\right).
    \label{eq:RewardFunction}
\end{equation}
Here, 
$\IndividualReward_\CurrentTime(\cdot)$ provides the number of orders fulfilled by a courier by taking the input action  
at time $\CurrentTime$, whereas 
$\ActionDelay_{\CurrentTime}(\cdot)$ represents the time required for a courier to deliver all the orders in its queue, including the ones associated with its current task 
and return to the warehouse. To ensure that the first term 
has more weight than the second one 
in the objective function, the multiplier of $\PParameter$ is introduced. This parameter serves as a constant which takes into consideration elements such as maximum allowed queue size and geographical area. Its computation involves determining the longest conceivable queue duration, encompassing travel duration between different points and from those points to the warehouse on the map. For instance, if we consider a maximum queue size of three, $\PParameter$ is determined by computing the three longest travel durations between different locations on the map. The incorporation of $\PParameter$ is thus aimed at giving precedence to the maximization of the total number of orders fulfilled in each time interval. Nevertheless, when two feasible actions serve an equal number of orders, the decision rests on the option that enables the courier to finish their tasks in the least amount of time. This emphasis ensures that couriers become available more swiftly to handle new groups of orders.

\subsection{Exogenous Information and Transition Function}
During each time step within the decision horizon, the system receives a collection of online orders which constitute the exogenous information. The orders arriving between time $\CurrentTime$ and $\CurrentTime + 1$ are denoted as $\ExogenousInformation_{\CurrentTime+1}$. Moreover, $\ExogenousInformation_0$ represents the orders which have accumulated during the time between the final time step in the previous day and the initial time step $\CurrentTime=0$, given a 24-hour planning horizon. Subsequently, the evolution of the system state from time $\CurrentTime$ to $\CurrentTime + 1$ is determined by the transition function that depends on the arrival of online orders and the decision tuple $\DecisionTuple_\CurrentTime \in \DecisionSet_\CurrentTime(\StateNote_\CurrentTime)$. By introducing the post-decision state \citep{powell2007approximate}, the state transition can be divided into two distinct parts. The post-decision state captures the system state immediately after a decision has been made but prior to the arrival of exogenous information in the subsequent time step. The initial transition \eqref{eq:FirstTrans} leads to the post-decision state via the action $\DecisionTuple_\CurrentTime$, and is denoted by $\PostDecisionState_\CurrentTime$. As unfulfilled orders exit the system at each time step, the post-decision state consists solely of information related to the couriers. The subsequent transition \eqref{eq:SecondTrans} occurs from the post-decision state to the next state, influenced by the arrival of exogenous information $\ExogenousInformation_{\CurrentTime+1}$:
\begin{subequations}
\label{m:Transitions}
    \begin{alignat}{2}
    & \PostDecisionState_\CurrentTime=\FirstTransition(\StateNote_\CurrentTime,\DecisionTuple_\CurrentTime)  \label{eq:FirstTrans} \\
    & \StateNote_{\CurrentTime+1}=\SecondTransition(\PostDecisionState_\CurrentTime,\ExogenousInformation_{\CurrentTime+1})  \label{eq:SecondTrans}
    \end{alignat}
\end{subequations}
Due to the assumption that unassigned orders exit the system at the end of each decision epoch, the state of orders at time $\CurrentTime + 1$ is defined as follows:
\begin{equation}
    \label{eq:OrderTransition}
    \OrdersState_{\CurrentTime+1} = \ExogenousInformation_{\CurrentTime + 1}
\end{equation}
Furthermore, the state of couriers at time $\CurrentTime + 1$ is described as follows:
\begin{equation}
    \label{eq:CourierTransition}
    \CouriersState_{\CurrentTime + 1} = \StateNote^{\CourierPost}_{\CurrentTime}
\end{equation}
such that $\StateNote^{\CourierPost}_{\CurrentTime}$ denotes the state of all couriers after taking the actions $\DecisionTuple_\CurrentTime$ and simulating their movements forward in time by one period, prior to the arrival of new exogenous information. For a courier, their state remains unchanged in the next time step if they are not on their shift or if they are at the warehouse without any assigned orders. Yet, if a courier is at the depot and receives new orders, they are sent out with their state updated to show the estimated return time. Similarly, if a courier is already out delivering and receives new online orders, their queue adapts while they continue their ongoing delivery route.
\subsection{Optimal Policy}
The objective in our order dispatching problem is to maximize the expected number of online orders served throughout the operation horizon:
\begin{equation}\label{eq:NeurADPObj}
    \max_{\pi \in \Pi}\mathbb{E}_{\ExogenousInformation=(\ExogenousInformation_0, \ldots, \ExogenousInformation_{T})} \left[ \sum_{\CurrentTime \in \Time} \Reward_\CurrentTime \left( \DecisionSet_\CurrentTime^\pi(\StateNote_\CurrentTime^\pi(\ExogenousInformation)) \right) \big|\StateNote_0 \right].
\end{equation}
Through the solution of Equation~\eqref{eq:NeurADPObj}, we can identify a policy $\pi$ from a set of feasible policies $\Pi$ which maximizes the reward when its recommended actions $\DecisionSet_\CurrentTime^\pi(\StateNote_\CurrentTime)$ are sequentially implemented at realized states. The realized states are defined as follows:
\begin{subequations}
\label{m:Transitions2}
    \begin{alignat}{2}
    & \StateNote_0^\pi(\ExogenousInformation)=\StateNote_{0}  \label{eq:OptOne} \\
    & \StateNote_{\CurrentTime+1}^\pi(\ExogenousInformation)=\SecondTransition(\FirstTransition(\StateNote_\CurrentTime^\pi(\ExogenousInformation), \DecisionSet_\CurrentTime^\pi(\StateNote_\CurrentTime^\pi(\ExogenousInformation))),\ExogenousInformation_{\CurrentTime+1}),~\CurrentTime=0,\ldots,T-1  \label{eq:OptTwo}
    \end{alignat}
\end{subequations}
such that $\StateNote_0$ corresponds to the initialized couriers, as well as the orders which have accumulated during the time interval from the final time step in the previous day up to the starting point at $t=0$ within the initial state of the decision horizon, provided a 24-hour decision horizon. The future reward is determined by taking the expectation with respect to the stochastic process described by $\ExogenousInformation$. The actions and states encountered during each decision epoch depend solely on the revealed random variables up to that point, and accordingly the overall reward relies on the realization of the complete vector $\ExogenousInformation$. By solving the Bellman optimality equations, the optimal values $V_\CurrentTime(\StateNote_\CurrentTime)$ at each state $\StateNote_\CurrentTime$ can be calculated as
\begin{equation}
V_\CurrentTime(\StateNote_\CurrentTime)=\max_{\DecisionTuple_\CurrentTime \in \DecisionSet_\CurrentTime(\StateNote_\CurrentTime)} \big\lbrace \Reward_\CurrentTime(\DecisionTuple_\CurrentTime) + \mathbb{E}_{\ExogenousInformation_{\CurrentTime+1}}[V_{\CurrentTime + 1}(\StateNote_{\CurrentTime + 1})] \big\rbrace
\label{eq:Bellman}
\end{equation}
where $\StateNote_{\CurrentTime + 1} = \SecondTransition(\FirstTransition(\StateNote_\CurrentTime,\DecisionTuple_\CurrentTime),\ExogenousInformation_{\CurrentTime+1})$. To compute the value function $V_t(\cdot)$, a backward induction procedure can be employed, which involves working backward in time from the final epoch $T$~\citep{powell2007approximate}. This procedure considers the rewards associated with taking the optimal actions and the probabilities of transitioning between states. The recursive process continues until reaching the first stage, $\CurrentTime=0$. However, this approach becomes impractical even for small instances due to the requirement of enumerating all possible outcomes and actions. Accordingly, ADP-based methods can be used to solve such problems. 

\section{Solution Methodology} \label{solutionmethodology}
In this section, we first describe an ADP approach for the ODP, highlighting its handling of the curses of dimensionality, VFA and updating methods. 
Subsequently, we introduce NeurADP as the more suitable method for our order dispatching problem and discuss its main distinctions from the ADP approach. Lastly, we provide our adaptation of the NeurADP algorithm for our problem setting.

\subsection{ADP and VFA for Ultra-fast ODP}
\label{adpmethodology}
As outlined in Section~\ref{problemformulation}, feasible decisions in the ODP involve not only the acceptance and rejection of incoming orders, but also the assignment of accepted order batches to courier queues and the determination of their respective routes. 
Recall the feasible set of decisions denoted by $\DecisionSet_\CurrentTime(\StateNote_\CurrentTime)$.
The definition of this set ensures that couriers are on their shifts, have adequate capacity for assigned orders, have viable routes for timely delivery, including any previously assigned orders, and can return to the depot before their shift concludes. Given this, the optimal policy for the ODP may be obtained using the Bellman optimality equations defined in Equation~\eqref{eq:Bellman}. However, computing $V_\CurrentTime(\StateNote_\CurrentTime)$ exactly proves intractable for complex large-scale problems such as the ODP due to what~\citet{powell2007approximate} classifies as the ``three curses of dimensionality'', referring to the challenges of managing the state, action, and outcome space. More specifically, solving the Bellman optimality equation for a state $\StateNote_\CurrentTime$ requires computing the anticipated downstream reward. This involves multiplying the value of each possible outcome $\StateNote_{\CurrentTime + 1}$ by the probability determined by the exogenous information $\ExogenousInformation_{\CurrentTime+1}$. However, for large-scale problems such as the ODP, the outcome space becomes excessively large, resulting in the first curse of dimensionality. To overcome this, ADP utilizes the concept of post-decision states, dividing the dynamic programming equation into two parts:
\begin{subequations}
\label{m:BellmanBreakup}
    \begin{align}
    V_\CurrentTime(\StateNote_\CurrentTime)\ &=\max_{\DecisionTuple_\CurrentTime \in \DecisionSet_\CurrentTime(\StateNote_\CurrentTime)} \{ \Reward_\CurrentTime(\DecisionTuple_\CurrentTime) + V_\CurrentTime^\PostDecision(\StateNote_\CurrentTime^\CourierPost) \} \label{eq:BellOne} \\
    V_\CurrentTime^\PostDecision(\StateNote_\CurrentTime^\CourierPost)\  & =\mathbb{E}_{\ExogenousInformation_{\CurrentTime+1}}\big[V_{\CurrentTime + 1}(\StateNote_{\CurrentTime + 1})\big|\StateNote_\CurrentTime^\CourierPost\big]  \label{eq:BellTwo}
    \end{align}
\end{subequations}
To avoid enumerating the entire outcome space and evaluating future values, Equation~\eqref{eq:BellOne} establishes a deterministic optimality equation based on the post-decision state. Hence, it eliminates the need for such exhaustive computations. Equation~\eqref{eq:BellTwo} expresses the post-decision state value function as the expected value of downstream rewards, where $\StateNote_{\CurrentTime + 1}=\SecondTransition(\StateNote_\CurrentTime^\CourierPost,\ExogenousInformation_{\CurrentTime+1})$. However, the computational challenge arises from the high-dimensional state space, making it difficult to compute value functions for all feasible post-decision states. In the ODP, the post-decision state is influenced by various factors related to couriers and their potential states, including their current locations and shift-times, as well as their assigned orders and their respective deadlines. This complexity increase, known as the ``second curse of dimensionality'', is associated with the exponential increase in the state space. To address this, an approximation of the post-decision state value function, $V^{\PostDecision}_\CurrentTime (\StateNote_\CurrentTime^\CourierPost)$, can be used.

There exist several types of value function approximations~\citep{powell2010approximate}. 
One approach to VFA involves the utilization of lookup tables in conjunction with state-wise aggregation. More specifically, a unique entry is assigned to each state in the lookup table, applying varying levels of aggregation to state values and using weighted summations to improve the accuracy of VFAs obtained from the aggregation levels. 
The utilization of basis functions offers an alternative approach for performing VFAs. These functions serve the purpose of transforming the original state space into a typically lower-dimensional form, with the goal of capturing influential state features which impact their values. 

Another method employed for VFA is the dual heuristic approach, which relies upon the concept of marginal values. More specifically, rather than exclusively focusing on the inherent value associated with occupying a particular state, this method places greater importance on assessing how the value function changes concerning that state's derivative, hence enabling the prioritization of the rate of change in value rather than the absolute value itself. This often leads to more efficient problem-solving across a wide array of practical applications. Furthermore, this approach is particularly valuable for resource allocation problems, such as the ODP, where vector-valued decision problems, namely \eqref{eq:BellOne} in this framework, are typically addressed, e.g., using linear, nonlinear, or integer programming. The dual heuristic approach has found widespread application in the transportation domain, for problems such as the ride-pool matching problem \citep{yu2019integrated}, taxi-on-demand \citep{al2020approximate,simao2009approximate}, and crowd-shipping \citep{mousavi2021approximate}. 
However, we observe that this approach would have some important drawbacks when applied to the ODP. Since this observation motivates our proposal of instead employing NeurADP for the ODP, we next briefly explain the dual heuristic approach.

In its common practice, we would define the linear approximation of the courier-based post-decision value function using the linear decomposition of the function $\bar{V}^\PostDecision_\CurrentTime (\StateNote_\CurrentTime^\CourierPost)$. This function incorporates the courier vector attributes in the post-decision state, specifically $\{{\StateNote_{\CurrentTime \IndividualCourier}^\CourierPost}\}_{\IndividualCourier \in \CouriersSet}$, and is formally defined as follows:
\begin{equation}
    \bar{V}_\CurrentTime^\PostDecision (\StateNote_\CurrentTime^\CourierPost) := \sum_{\IndividualCourier \in \CouriersSet} \EstimatedV_{\CurrentTime \IndividualCourier}^\PostDecision \StateNote_{\CurrentTime \IndividualCourier}^\CourierPost \label{eq:ADPLinApprox}
\end{equation}
wherein $\EstimatedV_{\CurrentTime \IndividualCourier}^\PostDecision$ represents the expected down-stream reward associated with a courier being in the post-decision state of $\IndividualCourier$ at time $\CurrentTime$. This representation provides a considerable computational benefit, since rather than computing value functions for each possible post-decision state at the current time $\CurrentTime$, we need only $\big|\CouriersSet\big|$ variables, which are denoted as $\EstimatedV_{\CurrentTime \IndividualCourier}^\PostDecision$. In the subsequent ADP algorithm, where $\SingleDay$ represents the iteration number and is used to index all the relevant components, the update for $\EstimatedV_{\CurrentTime \IndividualCourier}^\PostDecision$ is defined as follows:
\begin{equation}
    \EstimatedV_{{\CurrentTime}\IndividualCourier}^{\PostDecision,{\SingleDay}}=(1-\StepSize^\SingleDay)~\EstimatedV_{{\CurrentTime}\IndividualCourier}^{\PostDecision,{\SingleDay-1}} + \StepSize^\SingleDay~\MarginalValue_{{\CurrentTime}\IndividualCourier}^{\PostDecision,\SingleDay}
    \label{eq:UpdateEquation}
\end{equation}
such that $\StepSize^\SingleDay$ represents the step-size at iteration $\SingleDay$ of the algorithm, $\EstimatedV_{{\CurrentTime}\IndividualCourier}^{\PostDecision,{\SingleDay-1}}$ signifies the current approximate value of $\EstimatedV_{\CurrentTime \IndividualCourier}^\PostDecision$, and $\MarginalValue_{{\CurrentTime}\IndividualCourier}^{\PostDecision,\SingleDay}$ represents the observed marginal values associated with having an additional courier of type $\IndividualCourier$ at time $\CurrentTime$. The partial derivative values $\MarginalValue_{{\CurrentTime}\IndividualCourier}^{\PostDecision,\SingleDay}$ may be obtained as the numerical derivative of the following MIP model:
\begin{subequations}
\begin{align}
\text{max} \quad & \Reward_{\CurrentTime+1}(\StateNote_{\CurrentTime+1}^\SingleDay, \DecisionTuple_{\CurrentTime+1}) + \sum_{\IndividualCourier \in \CouriersSet} \EstimatedV_{{\CurrentTime+1},\IndividualCourier}^{\PostDecision,{\SingleDay-1}} \StateNote^{\CourierPost, \SingleDay}_{\CurrentTime+1,\IndividualCourier} \\
\text{s.t.} \quad & \DecisionTuple_{\CurrentTime+1} \in \DecisionSet^{\texttt{MIP}}_{\CurrentTime+1}(\StateNote^\SingleDay_{\CurrentTime+1})
\end{align}
\end{subequations}
where $\DecisionSet^{\texttt{MIP}}_{\CurrentTime}(\cdot)$ represents an MIP formulation of the ODP feasible decisions. Such an MIP model would necessitate introducing decision variables for order acceptance/rejection, assignment of order batches to couriers, as well as those to decide courier routes along with various sets of constraints to ensure their feasibility. As such, this would not be a computationally viable approach, in particular due to the need to solve this MIP to optimality a large number of times. Therefore, in the literature, the common approach has been to transform the MIP into its linear programming (LP) relaxation, e.g., $\DecisionTuple_{\CurrentTime+1} \in \DecisionSet^{\texttt{LP}}_{\CurrentTime+1}(\StateNote^\SingleDay_{\CurrentTime+1})$, and using LP duals; in our case this would be dual values associated with the constraints pertaining to courier flow conservation.


\subsection{Motivation for NeurADP}
\label{adpdrawbacks}
The ADP methodology described above is not suitable for our ODP due to several reasons.
Firstly, as noted for $\DecisionSet^{\texttt{MIP}}_{\CurrentTime+1}(\cdot)$, our problem setting involves a complex decision-making process which necessitates a more intricate one-to-many matching between couriers and batches of online orders, as opposed to a straightforward one-to-one courier-order matching. Additionally, complex routing decisions have to be made for each courier adhering to respective order deadlines, further complicating the decision space. On the other hand, due to poor LP relaxations, updating the VFA parameters with the dual values of the matching LP is not a preferable option either. Furthermore, to mitigate the curses of dimensionality, ADP commonly employs an aggregated attribute space. For the ODP, this approach involves consolidating courier-related attributes, such as their locations, rather than considering each courier's individual state. However, crafting these state attributes manually requires domain expertise, which is challenging in the context of the ODP. This primarily stems from the complexity of the ODP as it involves numerous couriers with varying shift times, different numbers of orders, each having unique deadlines, and traveling along distinct trajectories. While such a manual approach may be reasonable in transportation problems such as the taxi-on-demand problem with a single request per driver, it becomes impractical when dealing with couriers who have personalized shifts and multiple orders to manage. Finally, although linear and piece-wise linear VFAs offer simplicity in their integration into MIP models, this simplicity may diminish modeling accuracy and representational power. This becomes particularly evident in intricate, high-dimensional problem settings characterized by non-linear dynamics and complex attribute dependencies, as observed in the ODP, rendering them inferior options for our specific problem setting.

NeurADP, introduced by \citet{shah2020neural} to address the one-to-many case of the ride-pool matching problem, is an innovative ADP-based algorithm explicitly designed to overcome the limitations of traditional ADP methods when dealing with large-scale problems. While both NeurADP and ADP aim to solve sequential decision-making problems by approximating the value functions of post-decision states, they differ in their approach. As mentioned, ADP typically relies on linear or piece-wise linear VFAs. In contrast, NeurADP utilizes a non-linear neural network-based VFA. This allows for an automatic compact low-dimensional state-space representation without the need for domain expertise for state-space aggregation as a means of dealing with high state space dimensionality. The neural network-based non-linear value function is then innovatively integrated into the MIP-based framework through a \textit{two-step decomposition}:   (1) the set of feasible actions for each courier is enumerated, (2) a matching integer program (much simpler than the aforementioned MIP consisting of all the decisions of the ODP) is solved over all couriers, with the values associated with each action integrated into the integer programming (IP) model as constants. Furthermore, rather than using  LP-based duals to update these approximations as in ADP, NeurADP leverages DRL techniques for updating its value function approximations. More specifically, the gradients associated with the network parameters are computed and adjusted by minimizing the L2-norm between the current value function estimate and a one-step projection of the return derived from the Bellman equation. To enhance stability, NeurADP incorporates off-policy updates along with DRL techniques such as the implementation of a target network and Double Q-learning \citep{van2016deep}. These additions further refine the algorithm and contribute to its improved performance.
We explain these concepts in our ODP adaptation in more detail next.

\subsection{NeurADP Solution Methodology}
\label{neuradpsolution}
We next detail the NeurADP algorithm for our problem setting. We first describe the two-step decomposition enabled by NeurADP and explain identifying feasible courier-order matchings and the IP model for obtaining optimal matching. Then, the NeurADP-based VFA is explained, which is followed up by the description of the overall algorithm. Lastly, a brief discussion on the neural network architecture is provided.

\subsubsection{Two-step decomposition}
At every decision epoch $\CurrentTime$, given the state of the system $\StateNote_\CurrentTime=(\CouriersState_\CurrentTime,\OrdersState_\CurrentTime)$, the NeurADP solution methodology begins with enumerating feasible matchings between couriers and incoming batches of online orders. To evaluate the feasibility of batching a set of orders together or delivering a single order by itself, we take into consideration whether a batch can be delivered before each order's respective drop-off deadline. In other words, a batching is feasible if there exists a viable route for a courier to deliver each order in the batch, as well as the orders it is currently assigned to, prior to each order's respective deadline, beginning from the warehouse. With respect to the orders present in $\OrdersState_\CurrentTime$, we first define $\BatchingSet_{\CurrentTime}$ to represent the set of all order batchings with the minimum batch size of one and the maximum batch size of the available capacity of an empty courier's queue, denoted by $\MaxQueueSize$. It is important to note that while we denote the maximum queue size to be equal for all couriers in our model, this simplification is made for the sake of notational clarity, and it may be varied for each courier in practice. To assess whether it is possible to match a courier $\IndividualCourier$ to a new order batch $\FeasibleBatching$, we consider the following constraints: 
\begin{subequations}
\label{m:FeasibilityConsts1}
    \begin{alignat}{2}
\IndividualCourier_\CourierShiftStart & \leq \CurrentTime  \label{eq:StartConst} \\ \IndividualCourier_\CourierShiftStart +  \ShiftLength & > \CurrentTime  \label{eq:EndConst} \\
\big| \IndividualCourier_{\CourierOrders}  \big| + \big| \FeasibleBatching \big| & \leq \MaxQueueSize
    \label{eq:QueueConst}
    \end{alignat}
\end{subequations}
Here, the constraints presented in equations~\eqref{eq:StartConst} and~\eqref{eq:EndConst} ensure that the courier is actively working during the time when the batch of orders is assigned to them, with $\ShiftLength$ denoting the shift length of couriers, while constraint \eqref{eq:QueueConst} guarantees that adding the orders from the newly assigned order batch to the courier's existing queue of previously matched orders does not exceed the maximum allowed capacity. Assuming the courier is currently on duty and has adequate storage capacity for both their new and previously assigned orders, we next confirm the existence of a feasible sequence for these combined orders. This sequence must ensure that each order delivery is completed before the respective drop-off deadline, and that the courier is able to successfully complete all deliveries and return to the warehouse prior to the end of their shift. We establish the new queue of online orders for courier $\IndividualCourier$, encompassing both the previously assigned orders as well as those within batch $\FeasibleBatching$, by $\IndividualCourier_{\CourierOrders'}$, such that $\IndividualCourier_{\CourierOrders'} = \IndividualCourier_\CourierOrders \cup \FeasibleBatching$. We then define $\QueueCadrinality=|\IndividualCourier_{\CourierOrders'}|$, and introduce $\PermutationSet$ as the set of all permutations of $\IndividualCourier_{\CourierOrders'}$ where each $\SinglePermutation \in \PermutationSet$ represents a unique sequence for delivering the online orders in $\IndividualCourier_{\CourierOrders'}$. 
For instance, if $\IndividualCourier_{\CourierOrders'} = \{4,6,12\}$ and $\SinglePermutation = (12, 4, 6)$, then order 12 from $\IndividualCourier_{\CourierOrders'}$ is delivered first, followed by the order 4, and finally the order 6. Furthermore, given permutation $\SinglePermutation$, we let $\IndividualOrder_{\IndividualVehicleDropOff}^{\SinglePermutation(\SingleIndex)}$ indicate the earliest drop-off time of the order at the drop-off index $\SingleIndex$ in $\SinglePermutation$, and $\IndividualOrder_\OrderDestination^{\SinglePermutation(\SingleIndex)}$ to represent its drop-off destination. Note that, given $\SinglePermutation$, these earliest drop-off times can be calculated in a forward manner starting from $\SingleIndex=1$ to $\SingleIndex = \QueueCadrinality$, where $\IndividualOrder_{\IndividualVehicleDropOff}^{\SinglePermutation(1)} = \CurrentTime + \TravelTime_{\CurrentTime}(\IndividualOrder, \IndividualOrder_\OrderDestination^{\SinglePermutation(1)})$, $\IndividualOrder_{\IndividualVehicleDropOff}^{\SinglePermutation(2)} = \IndividualOrder_{\IndividualVehicleDropOff}^{\SinglePermutation(1)} + \TravelTime_{\CurrentTime}(\IndividualOrder_\OrderDestination^{\SinglePermutation(1)}, \IndividualOrder_\OrderDestination^{\SinglePermutation(2)})$ and so on. Our goal is to ensure the existence of a permutation $\SinglePermutation \in \PermutationSet$ of online orders which meets the following constraints:
\begin{subequations}
\label{m:FeasibilityConsts2}
    \begin{alignat}{2}
    & \IndividualOrder_\RequestDeadline^{\SinglePermutation(\SingleIndex)} \leq \IndividualOrder_\IndividualVehicleDropOff^{\SinglePermutation(\SingleIndex)} && \quad \quad  \quad \forall \SingleIndex \in \{1,\ldots,\QueueCadrinality\}
    \label{eq:DelayConst} \\ 
    & \IndividualCourier_\CourierShiftStart +  \ShiftLength > \IndividualOrder_\IndividualVehicleDropOff^{\SinglePermutation(\QueueCadrinality)} + \TravelTime_{\CurrentTime}( \IndividualOrder_\OrderDestination^{\SinglePermutation(\QueueCadrinality)},0)
    \label{eq:ShiftConst}
    \end{alignat}
\end{subequations}
Constraint \eqref{eq:DelayConst} ensures that the courier is able to deliver all of their assigned orders prior to their individual deadlines. 
Furthermore, constraint \eqref{eq:ShiftConst} guarantees that the courier is able to fulfill all deliveries and return the depot prior to the end of their shift. 
Then, the feasible set of matches at time $\CurrentTime$ between couriers and order batchings may be described as follows:
\begin{equation}
\FeasibilitySet_\CurrentTime = \{(\IndividualCourier,\FeasibleBatching) \in \CouriersState_\CurrentTime \times \BatchingSet_\CurrentTime : \eqref{eq:StartConst}-\eqref{eq:QueueConst},~\eqref{eq:DelayConst}-\eqref{eq:ShiftConst}\}
\label{eq:FeasibilitySet}
\end{equation}
We note that the feasible set of matches, $\FeasibilitySet_\CurrentTime$, can be enumerated efficiently for the ultra-fast ODP with the capacitated couriers. Whereas for higher dimensional or less restricted problem settings (e.g., see \citep{shah2020neural}), full enumeration might not be achievable, in which case a subset of feasible matches can be heuristically generated.

In the second step, the NeurADP algorithm builds the following matching IP model to determine the decisions for each courier:
\begin{subequations}
\label{m:LinearSystem}
    \begin{alignat}{2}
    \texttt{MatchingIP:} \quad \text{max} \ & \sum_{\IndividualCourier \in \CouriersState_{\CurrentTime}} \sum_{\MovementPoint \in \FeasibilitySet_\CurrentTime \cup \{\NullAction\}} \ImmediateReward_{\CurrentTime\IndividualCourier\MovementPoint} \cdot \MatchingVariable_{\CurrentTime \IndividualCourier \MovementPoint} + \texttt{score}_{\CurrentTime}(\IndividualCourier \leftrightarrow \MovementPoint)\cdot \MatchingVariable_{\CurrentTime \IndividualCourier \MovementPoint} \label{eq:ObjectiveValue} \\ 
    \text{s.t.} \ & \sum_{\FeasibleBatching \in \BatchingSet_\CurrentTime: (\IndividualCourier, \FeasibleBatching) \in \FeasibilitySet_\CurrentTime} \MatchingVariable_{\CurrentTime\IndividualCourier\FeasibleBatching} + \MatchingVariable_{\CurrentTime \IndividualCourier \NullAction} = 1 && \quad \quad \quad \forall \IndividualCourier \in \CouriersState_\CurrentTime
    \label{eq:VehicleFlow} \\
    & \sum_{\IndividualCourier \in \CouriersState_\CurrentTime}\sum_{\FeasibleBatching \in \BatchingSet_\CurrentTime: (\IndividualCourier, \FeasibleBatching) \in \FeasibilitySet_\CurrentTime;\IndividualOrder \in \FeasibleBatching} \MatchingVariable_{\CurrentTime\IndividualCourier\FeasibleBatching} \leq 1 && \quad \quad \quad \forall \IndividualOrder \in \OrdersState_\CurrentTime \label{eq:PassengerFlow} \\
    & \MatchingVariable_{\CurrentTime\IndividualCourier\MovementPoint} \in \{0,1\} && \hspace{-1.2cm}\forall \IndividualCourier \in \CouriersState_\CurrentTime, \MovementPoint \in \FeasibilitySet_\CurrentTime ~\cup \{\NullAction\} \label{eq:VariableDefinition}
    \end{alignat}
\end{subequations}
At each time step, there exists a default action for each courier denoted by $\NullAction$, indicating that they will not be assigned a new batch of orders.
The constraints represented by~\eqref{eq:VehicleFlow} guarantee that each individual courier $\IndividualCourier$ at time $\CurrentTime$ is assigned to exactly one feasible action. Similarly,~\eqref{eq:PassengerFlow} ensures that each individual order $\IndividualOrder$ is assigned to at most one courier. Furthermore,~\eqref{eq:VariableDefinition} ensures that all decision variables are binary. (Importantly, we note that the set of actions for feasible matchings, $\FeasibilitySet_\CurrentTime$, may be further reduced by considering only the permutation for each updated queue $\IndividualCourier_{\CourierOrders'}$ which minimizes the time taken to make all deliveries in the queue and return to the depot.) Finally, the objective in~\eqref{eq:ObjectiveValue} calculates the total reward over the immediate reward of assigning courier $\IndividualCourier$ to the feasible matching $\MovementPoint$, which is calculated as in Equation~\eqref{eq:RewardFunction} and denoted by $\ImmediateReward_{\CurrentTime\IndividualCourier\MovementPoint}$, as well as the downstream reward, denoted by $\texttt{score}_{\CurrentTime}(\IndividualCourier \leftrightarrow \MovementPoint)$, gained from matching a courier $\IndividualCourier$ to a feasible decision $\MovementPoint$ at time $\CurrentTime$.

\subsubsection{NeurADP-based VFA}
Next, given the \texttt{MatchingIP}, we detail how NeurADP approximates the value functions, linking it to the ADP content reviewed in Section \ref{adpmethodology}. 
To derive the dispatching policy, NeurADP looks to solve the Bellman optimality equations introduced in Equation~\eqref{eq:Bellman}. Similar to ADP, NeurADP utilizes the concept of post-decision states, and divides the dynamic programming equation into two parts, namely Equation~\eqref{eq:BellOne} and Equation~\eqref{eq:BellTwo}. Furthermore, it aims to approximate the post-decision state value function, denoted by $V_\CurrentTime^\PostDecision(\StateNote_\CurrentTime^\CourierPost)$, also by first performing a courier-based decomposition. More specifically, the value function of the courier-based post-decision state is decomposed into the individual couriers' value functions as follows:
\begin{equation}
\label{eq:CourierDecomposition} \bar{V}_\CurrentTime^{\PostDecision}(\StateNote_\CurrentTime^\CourierPost) \approx \sum_{\IndividualCourier \in \CouriersSet} \widetilde{V}_{\CurrentTime \IndividualCourier}^\PostDecision(\StateNote_{\CurrentTime \IndividualCourier}^\CourierPost).
\end{equation}
NeurADP then approximates the value functions of individual couriers. In doing so, the assumption is made that a courier's long-term reward is minimally affected by the actions of other couriers in the current decision epoch. This assumption, rooted in the idea that long-term rewards primarily stem from the interaction between delivery routes, enables modeling that focuses on the pre-decision state of other couriers. Thus, it simplifies and expedites the optimization of delivery routes and schedules by reducing computational complexity and not heavily weighing the numerous possible actions of other couriers. The approximation of individual courier value functions can be succinctly expressed as:
\begin{equation}
\label{eq:CourierApproximation} \widetilde{V}_{\CurrentTime \IndividualCourier}^{\PostDecision}(\StateNote_{\CurrentTime \IndividualCourier}^\CourierPost) \approx \hat{V}_{\CurrentTime \IndividualCourier}^\PostDecision(\StateNote_{\CurrentTime \IndividualCourier}^\CourierPost,\{\StateNote_{\CurrentTime \IndividualCourier'}\}_{\IndividualCourier' \neq \IndividualCourier}).
\end{equation}
Here, the approximated post-decision state value function for courier $\IndividualCourier$, denoted by $\bar{V}_{\CurrentTime \IndividualCourier}^\PostDecision$, accepts input data about its own post-decision state, $\StateNote_{\CurrentTime \IndividualCourier}^\CourierPost$, as well as auxiliary pre-decision state information from the other couriers, $\{\StateNote_{\CurrentTime \IndividualCourier'}\}_{\IndividualCourier' \neq \IndividualCourier}$. This auxiliary data provides context about the environment in which courier 
$\IndividualCourier$ operates before taking an action, including the number of couriers on break, those at the warehouse, the average occupied capacity of other couriers, and the volume of incoming orders at that time step. Incorporating this additional information allows for a more accurate evaluation of the post-decision state value, considering that order acceptance is influenced not just by a courier's state, but also by the competitive nature of their operational environment. 
The overarching value function may thus be written as follows:
\begin{equation}
\bar{V}_\CurrentTime^{\PostDecision}(\StateNote_\CurrentTime^\CourierPost) \approx \sum_{\IndividualCourier \in \CouriersSet}\hat{V}_{\CurrentTime \IndividualCourier}^\PostDecision(\StateNote_{\CurrentTime \IndividualCourier}^\CourierPost,\{\StateNote_{\CurrentTime \IndividualCourier'}\}_{\IndividualCourier' \neq \IndividualCourier}).
\label{eq:CourierFinal}
\end{equation}

The individual value functions, $\hat{V}_{\CurrentTime \IndividualCourier}^\PostDecision(\cdot)$,  
are linearly integrated into the overall value function within the \texttt{MatchingIP}. That is, $\texttt{score}_{\CurrentTime}(\IndividualCourier \leftrightarrow \MovementPoint)$ terms in the \texttt{MatchingIP} objective, which reflect the estimated long-term value of assigning a specific courier ($\IndividualCourier$) to a particular feasible matching ($\MovementPoint$), are derived from $\hat{V}_{\CurrentTime \IndividualCourier}^{\PostDecision}(\cdot)$. Specifically, for each $(\IndividualCourier, \MovementPoint)$ pair, first the post-decision state, $\StateNote_{\CurrentTime \IndividualCourier}^\CourierPost$, is determined and then the corresponding value from the approximated value function, $\hat{V}_{\CurrentTime \IndividualCourier}^{\PostDecision}(\cdot)$, is obtained from the trained neural network to calculate the $\texttt{score}_{\CurrentTime}(\IndividualCourier \leftrightarrow \MovementPoint)$ values.
This approach reduces the evaluations of the non-linear value function from an exponential to a linear scale with respect to the number of couriers. 

To update the individual courier-based VFAs, NeurADP explicitly calculates the gradients associated with each parameter using standard symbolic differentiation libraries. It then adjusts these parameters to minimize the L2 distance between the one-step return estimate of the Bellman equation and the current value function estimate. In order to tackle stability and scalability issues, which are particularly vital in neural network value function learning, NeurADP employs a combination of methodological and practical strategies. These include the use of off-policy updates to stabilize Bellman updates 
and addressing data scarcity by directly storing sets of feasible actions. NeurADP additionally utilizes a singular neural network for individual courier value functions and employs prioritized experience replay to reuse experience efficiently. Moreover, practical simplifications such as utilizing low-dimensional embeddings for discrete locations and introducing Gaussian noise for exploration during training are strategically implemented. This ensures that learning remains manageable and is precisely tailored to the complexities and subtleties of the underlying problem space.

\subsubsection{Overall algorithm} 
\begin{figure}[!ht]
\begin{center}
\includegraphics[width=0.909\textwidth]{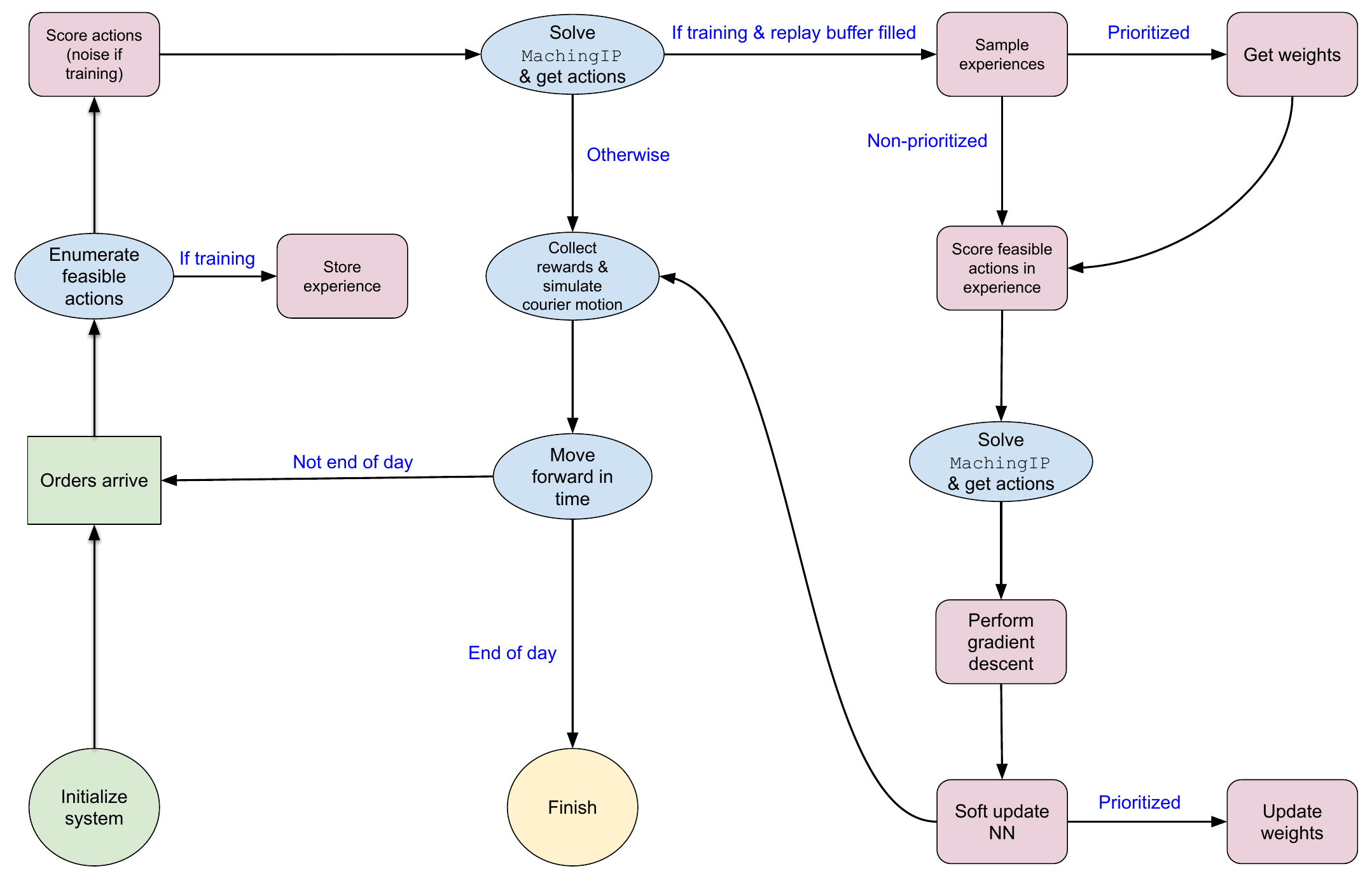}
\end{center}
\caption{NeurADP algorithm flowchart.}
\label{fig:NeurADPflowchart}
\end{figure}

The overall NeurADP algorithm for the ODP is presented in Figure~\ref{fig:NeurADPflowchart} in the form of a flow-chart. Initially, the system, along with prediction and target neural networks (used as value functions for courier post-decision states), and couriers with their shifts are initialized. Orders stochastically arrive and, depending on courier availability, feasible actions between couriers and orders are enumerated. During training, these actions are stored as future training experiences. The prediction neural network scores each action based on immediate rewards and the value of the resultant post-decision state, with Gaussian noise added for exploratory purposes during training \citep{plappert2017parameter}. 

The \texttt{MatchingIP} defined in \eqref{eq:ObjectiveValue}-\eqref{eq:VariableDefinition} is utilized to determine optimal actions, maximizing immediate and anticipated downstream rewards. If training, the replay buffer is checked for sufficient experiences to begin sampling and value function training. When utilizing a prioritized replay buffer, the associated weights with each experience are retrieved. The target neural network scores feasible actions from the experience, and the \texttt{MatchingIP} chooses the best actions. The prediction network updates the value of each post-decision state at time $\CurrentTime$ through gradient descent, using the \texttt{MatchingIP}-selected best action value at time $\CurrentTime + 1$. If a prioritized replay is in use, experience weights are updated. After sampling experiences, the algorithm collates rewards from the current iteration, simulates courier movements, and advances in time. If the subsequent time step marks the day's end, the iteration concludes; otherwise, it progresses to the new time step. 

\subsubsection{Neural network architecture} The architecture of the underlying neural network value function begins with an embedding layer which takes as input the current location of the courier as well as the destination locations of its matched orders. From here, these embedded location representations, complemented by their associated delays, are inputted into an LSTM layer. The output is then combined with additional pertinent auxiliary information and proceeds through several dense layers, ultimately yielding a single value. Furthermore, parameter tuning is undertaken throughout the network’s architecture, encompassing modifications to embedding sizes and variations in the number of dense layers.

\section{Experimental Setup} \label{experimentalsetup}
In this section, we describe the datasets and the benchmark policies employed in our numerical study.
All the experiments are conducted using Python 3.6.13 on Google Cloud servers and we use IBM ILOG CPLEX Optimization Studio version 12.10.0 to solve the IP models. 
\subsection{Datasets} 
In our numerical study, we examine four distinct geographical-based datasets containing delivery information of online orders: the Brooklyn, Chicago, Bangalore, and Iowa datasets. The Brooklyn and Chicago datasets~\citep{barkingdata2022} encompass delivery data for online DoorDash requests within their respective urban cities, while the Bangalore dataset~\citep{poddar2019} incorporates order requests from restaurants in Bangalore, India. Finally, the Iowa dataset, which was introduced by \citet{ulmer2021restaurant}, is comprised of destination locations for meal deliveries within Iowa City. Each dataset includes latitude-longitude points for drop-off locations associated with real-world order requests from which we extract a frequency distribution of the popularity of each delivery location within the city from which we sample. The datasets for Brooklyn, Iowa, Chicago, and Bangalore consist of 988, 500, 117, and 77 unique destination locations, respectively. A warehouse is located at the centre of each distribution of order destinations and the distance between each pair of locations is calculated via the haversine formula. To account for real-world travel time variations influenced by travel direction and accurately depict the unevenness in travel between two points, randomness is incorporated into the travel time values. This involves introducing an additional noise of up to 10\% to each value. 

\begin{figure}[!ht]
\begin{center}
\includegraphics[width=0.659\textwidth]{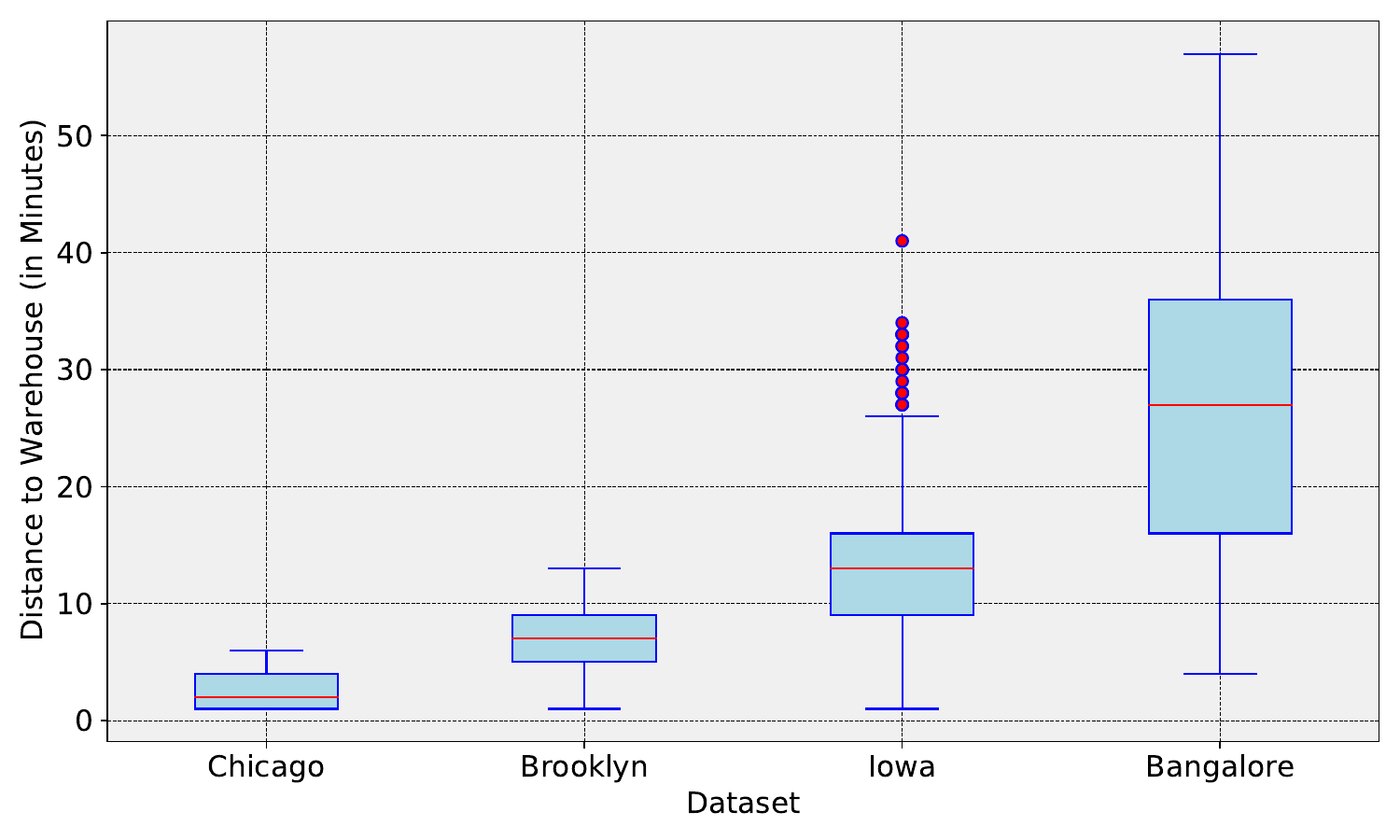}
\end{center}
\caption{Distribution of drop-off points for each dataset.}
\label{fig:data_boxplot}
\end{figure}

We consider a 24-hour problem horizon which is broken into 5-minute decision epochs (i.e., $\TimeIntervals = 5$ minutes). For all the datasets, the number of requests which arrive between each decision epoch is sampled based upon a distribution of real-world order requests, as presented in \citet{kavuk2022order}, with a 1-order request standard deviation band. Figure~\ref{fig:data_boxplot} shows a series of box plots that illustrate the distribution of distances between drop-off locations and the warehouse in each dataset. The calculation of the multiplier term $\PParameter$ in the reward function~\eqref{eq:RewardFunction} depends on the courier capacity parameter and the dataset as described in Section~\ref{sec:DecisionVariables}. For instance, for the Brooklyn dataset, for the capacity of 1, 2, 3 and 4, the $\PParameter$ values are 50, 73, 96, and 119, respectively. Similarly, for the courier capacity value of 3, the $\PParameter$ values are 467, 43, and 260 for Bangalore, Chicago and Iowa datasets, respectively.

\begin{figure}[!ht]
    \centering
    \subfloat[Location points\label{subfig:location_points}\centering]{{\includegraphics[width=0.405\textwidth]{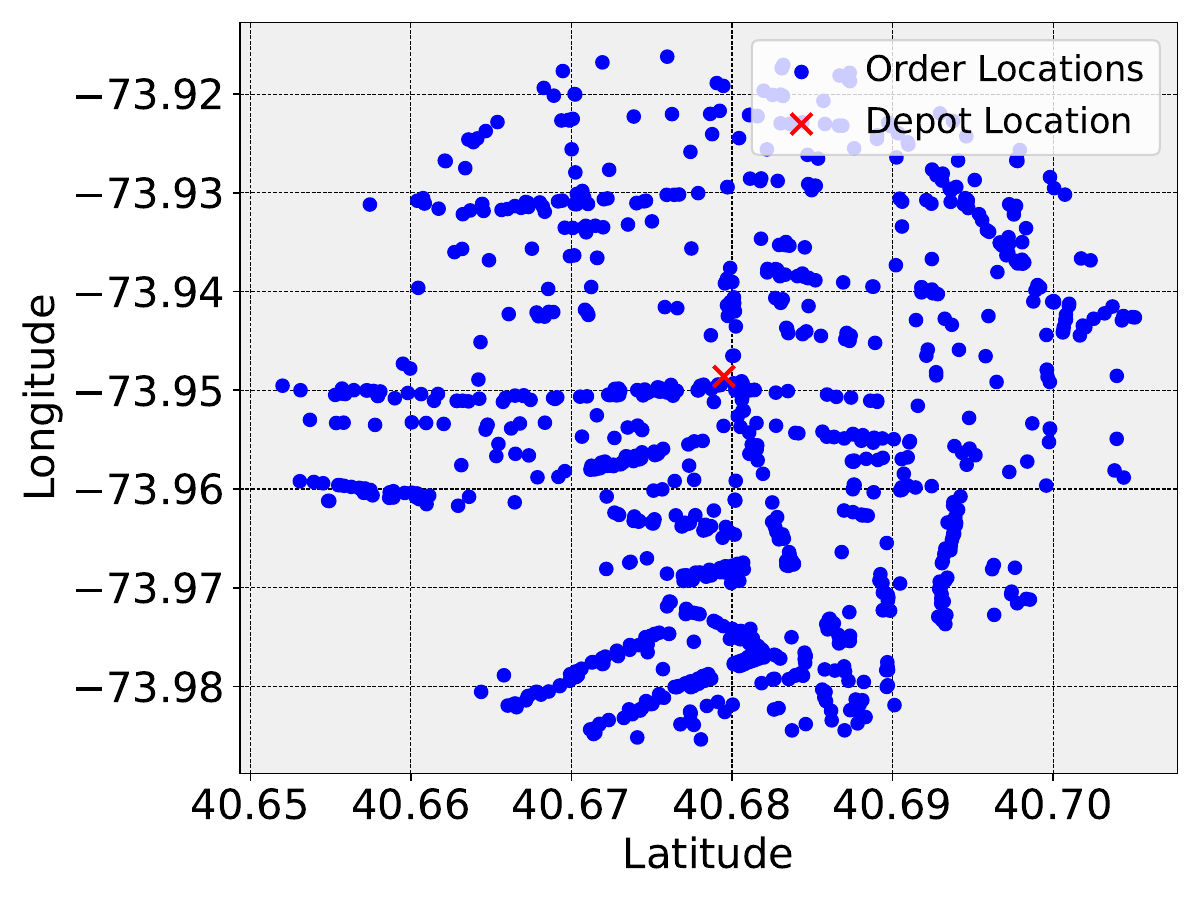} }}
    \subfloat[Order request distribution\label{subfig:order_distribution}\centering]{{\includegraphics[width=0.405\textwidth]{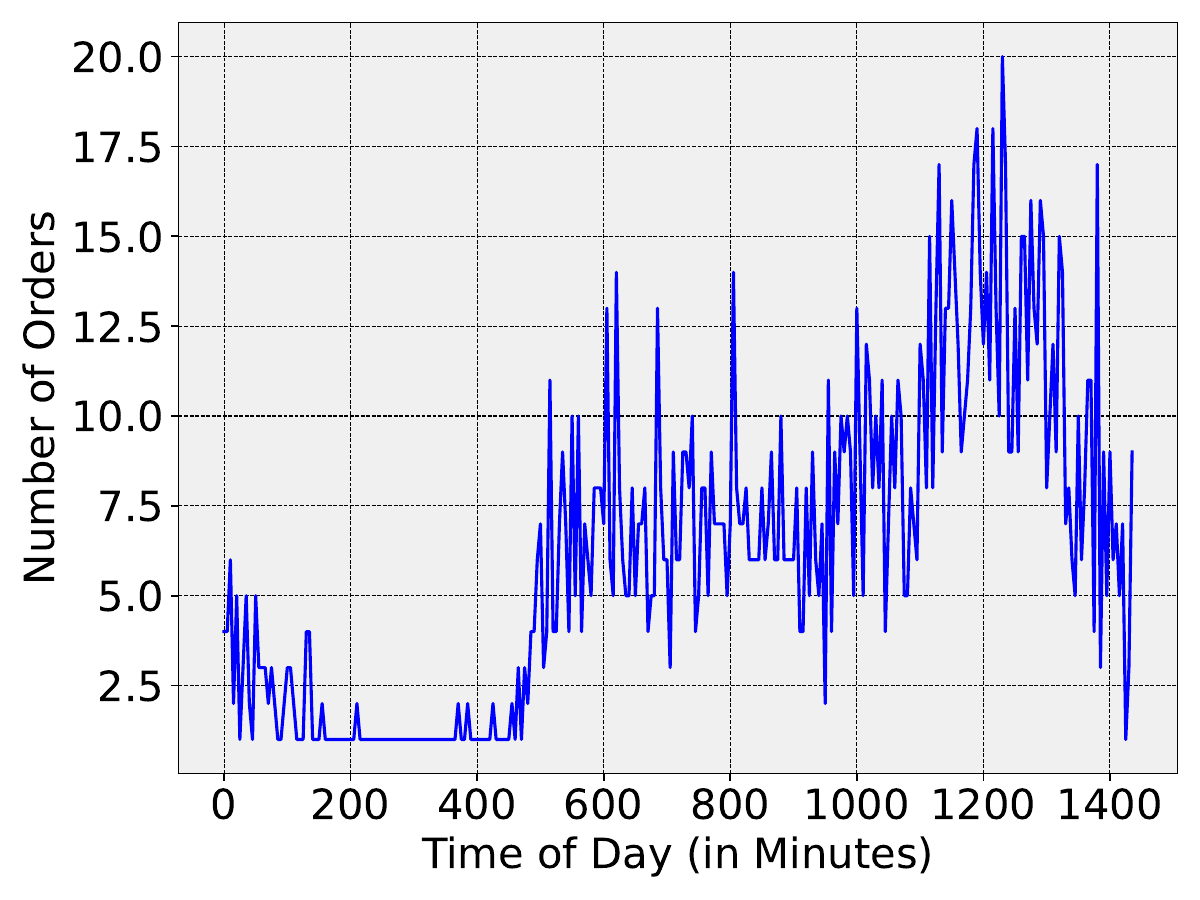} }}%
    \caption{Brooklyn dataset specifications.}
    \label{fig:data_summary}
\end{figure}

Figure~\ref{fig:data_summary} displays the geographic distribution of delivery locations in the Brooklyn dataset (as a representative dataset) along with the distribution of order arrivals throughout the day. The schedules of couriers are manually planned in advance to accommodate the anticipated fluctuations in order volume throughout the day. This entails scheduling fewer shifts during expected periods of low demand and scheduling more shifts during peak hours. More specifically, the quantity of couriers on duty closely mirrors the pattern of order arrivals depicted in Figure~\ref{subfig:order_distribution}. Between the hours of 3 AM and 6 AM, when the order volume is at its lowest, the courier staffing reaches its minimum. Conversely, the staffing level reaches its peak during the rush-hour period of 7 PM to 9 PM.

\subsection{Benchmark Policies}
We consider a family of benchmark policies in our comparative analysis with \NeurADP, namely, a set of myopic policies, which we denote by \Myopic, as well as a group of DRL policies, labeled as \DRL. Generally, myopic policies encompass greedy strategies which prioritize immediate rewards obtained from actions taken in the present time step, while disregarding any future consequences of these decisions. These policies facilitate a streamlined decision-making process, offering benefits when addressing complex and dynamic problems where it may not be feasible to calculate a globally optimal policy, as seen in the case of our order dispatching problem. However, as a result of their greedy nature, myopic policies may not always produce optimal long-term policies, and are thus often used as baselines. The employed \Myopic~policies follow a similar pattern and are implemented as follows. Available couriers are sorted so as to prioritize varying system dynamics, such as distance to the warehouse and available capacity. Each courier is then examined individually, with the aim of identifying the action between the courier and the set of incoming orders which maximizes the courier's order fulfillment while minimizing delivery time. For instance, in the case where couriers are sorted based upon proximity to the warehouse, we begin by evaluating feasible actions for the nearest courier. Among the available actions, we select the one which maximizes the number of orders which are matched, while minimizing the delivery time for both newly assigned and previously assigned orders. In cases where multiple actions serve an equal number of orders, the preference is given to the action which allows the courier to complete their assigned deliveries more quickly. Once an action is determined for a courier, we finalize the matching and move on to the next available courier in the sorted queue.

We examine a collection of \DRL~policies the ultra-fast ODP, which are derived from those introduced in \citep{kavuk2022order}. At each time step, we employ a trained Double Deep-Q Network (DDQN) to make accept-reject determinations for incoming online orders. The problem and decision dynamics are outlined as follows. The orders are sorted based upon their estimated delivery durations (direct delivery time). We evaluate each order individually and utilize the DDQN network to determine whether to accept or reject the order. If the decision is to accept, we employ a straightforward heuristic to match the order with an available courier, taking into account varying system dynamics such as distance to the warehouse and available capacity. A reward of 1 is accrued if an order is accepted and successfully matched with an available courier, and 0 reward is accrued otherwise. In the event that an order is accepted but cannot be fulfilled due to courier unavailability caused by capacity or time constraints, the order is disregarded, resulting in no reward. Once a decision is made and the order is either matched or ignored, we proceed to the next order in the sorted queue. The DDQN network, responsible for learning the accept-reject actions, is a feed-forward neural network and it is trained with the same dataset used for \NeurADP. It receives inputs regarding the state of the couriers, including their distances to the warehouse, their currently assigned orders, and the number of couriers on not on shift or at the warehouse. Additionally, information about the specific order under consideration is incorporated, encompassing details such as the destination location and delivery deadline. The network comprises four hidden layers with 32, 64, 64, and 32 neurons, respectively, and an output layer with two neurons representing the accept and reject decisions.

\section{Results}\label{results}
In this section, we present results from our detailed numerical study for the ultra-fast ODP. We primarily consider the average number of orders fulfilled within the 24-hour decision horizon while comparing the performance of the \NeurADP~policy against the two classes of benchmark policies. More specifically, we assess each policy over 20 days of testing data and subsequently compute the average total number of orders that each policy has encountered and fulfilled across these test days. Orders are generated via sampling at every decision epoch, and, according to a given policy, a feasible decision is made to batch them together and match them with couriers. Upon matching, the couriers' assigned orders are updated, and a simulation of the couriers' movement toward their next destination takes place. Below, we first introduce two novel artificial bounds on the achievable performance level that enhance our understanding of the quality of our solutions. Then, we delve into the analysis of our overall findings.

\subsection{Artificial Bounds on Achievable Performance Level}
Considering that problems such as the ODP often involve large, complex and highly dynamic systems, deriving optimal solutions for such problems is often infeasible and methodologies for deriving theoretical bounds are usually unsatisfactory \citep{bertsekas2005dynamic}. As such, practitioners typically evaluate policy performance based on the total potential reward available throughout the problem horizon. In the case of the ODP, this corresponds to the total number of online orders received in the system. However, depending upon the problem context and parameters, this overarching upper limit may not be realistic or reasonable to achieve, as it may be infeasible to come close to the maximum potential reward. For instance, in the ODP, though the system may receive a large number of orders throughout the day, provided the capacity constraints of couriers, order deadlines, number of couriers, and courier shift times, it may not be feasible to serve even half of those orders. As such, we introduce two novel benchmark ceilings. 

For our first benchmark ceiling, we consider a scenario where orders are served immediately upon being matched with couriers, which we define as \Direct. Specifically, we maintain the same constraints based on courier capacities and order deadlines. However, once a batch of orders is matched with the available couriers, we assume that all orders are delivered by the start of the subsequent decision epoch. This implies that all available couriers, as long as they are on their shift, will be ready in the following time step with an empty queue. This allows us to establish a lower bound, relative to the total number of orders observed, on how well our \NeurADP~and benchmark policies could have performed, even in this unrealistic setting. The second benchmark ceiling involves applying the NeurADP framework to each individual day's worth of orders in our testing dataset. Instead of training our \NeurADP~policy on a separate training dataset, deriving a policy that maximizes the expected reward, and subsequently testing it on a separate testing dataset, we train the \NeurADP~policy on the same testing dataset that we ultimately evaluate it on. This results in an alternative \NeurADP~policy, which we refer to as \NeurADPFixed. This policy aims to derive the best policy based on a deterministic set of orders for a specific day. This approach provides us with additional insights into the performance of alternative policies for ultra-fast ODP.

\subsection{Order Dispatching Performance}
We evaluate the results from our experiments with respect to four primary inputs: the number of couriers, delay time, courier capacity, and geographical location. The number of couriers is obtained by accounting for all the couriers working within a 24-hour period, while the delay time represents the maximum duration of time a courier has from an order's entry into the system to their drop-off. This duration is used to determine the order deadline, calculated using Equation~\eqref{eq:OrdDeadline}. The courier capacity specifies the maximum number of orders a courier can carry simultaneously, while the geographic location pertains to the spatial dataset and the distribution of requests based on geography. In our \textit{baseline configuration}, we utilize the Brooklyn dataset and set the parameters to include 15 couriers, a maximum allowable delay time of 10 minutes, and a maximum courier capacity of 3 orders. Due to the extensive computational time required to derive the \NeurADPFixed~ceiling values for each experiment, we consider it exclusively for the experiments related to our baseline configuration and the varying number of couriers and utilize the \Direct~ceiling for the remaining experiments. We begin by examining the baseline configuration to identify the most suitable benchmark policies from the \Myopic~and \DRL~policy classes, which are later used in the comparative analysis with the \NeurADP~policy.

\subsubsection{Baseline Configuration}
Our primary benchmark policies exhibit several variations in the matching process between orders and couriers. In the \DRL~policy, we have the flexibility to match accepted orders with couriers based on either their distance or current queue capacity. In terms of distance, the order can be assigned to the closest available courier or the farthest one. Regarding capacity, we can allocate the order to the courier with the least occupied queue or the most occupied queue. Similarly, for the \Myopic~policy, the matching of a batch of orders can be determined by considering both the proximity of distance and the capacity of the courier queue. We consider four variations of each benchmark policy. ``DC'' represents the utilization of distance in the matching process, where the closest courier is chosen. Conversely, ``DF'' signifies that the farthest courier is selected based on distance. Moreover, ``CE'' denotes the utilization of capacity for matching, with the emptiest courier being selected, while ``CF'' indicates the selection of the fullest courier based on capacity.

Table~\ref{table:base_table} provides the outcomes of the considered policy variants for the baseline configuration of our experimental setup. The table includes the \NeurADPFixed~value, which indicates the average number of orders served using the fixed ceiling, as well as the percentage of orders fulfilled by each policy, denoted as ``\% Filled''. More specifically, we calculate the average number of orders served by each policy and divide it by the fixed ceiling value. This result is then multiplied by 100, and a standard deviation is provided for each policy. Furthermore, the percentage increase of the \NeurADP~policy compared to the other benchmark policies is included in the final right-most column labeled ``\% Incr. NeurADP''. This metric is calculated by subtracting the average number of orders fulfilled by the benchmark policies from the average number of orders fulfilled by the \NeurADP~policy, then dividing the result by the ``Fixed Ceiling'' value, multiplied by 100.

\setlength{\tabcolsep}{4.5pt}
\renewcommand{\arraystretch}{1.15}
\begin{table}[!ht]
\centering
\caption{Performance of different policies for ultra-fast ODP for \textit{baseline configuration} (avg. number of orders fulfilled over 20-day test window is reported for the Fixed Ceiling; performance of other policies are w.r.t. Fixed Ceiling, provided as mean$\pm$stdev)}.
\label{table:base_table}
\resizebox{0.575\textwidth}{!}{
\begin{tabular}{P{0.15\textwidth}L{0.10\textwidth}L{0.20\textwidth}C{0.18\textwidth}}
\toprule
\textit{\textbf{Policy}} & \textit{\textbf{Fixed Ceiling}} & \textit{\textbf{\% Filled}} & \textit{\textbf{\% Incr. NeurADP}}\\
\midrule
\textit{\textbf{NeurADP}} & 955.00 & 97.96 $\pm$ 1.97 & -\\
\midrule
\textit{\textbf{DRL-DC}} & - & 86.04 $\pm$ 1.53 & +11.92\\
\textit{\textbf{DRL-DF}} & - & 83.65 $\pm$ 1.74 & +14.31\\
\textit{\textbf{DRL-CE}} & - & 81.04 $\pm$ 1.75 & +16.92\\
\textit{\textbf{DRL-CF}} & - & 85.81 $\pm$ 1.64 & +12.15\\
\midrule
\textit{\textbf{Myopic-DC}} & - & 91.27 $\pm$ 1.64 & +6.69\\
\textit{\textbf{Myopic-DF}} & - & 87.79 $\pm$ 2.15 & +10.17\\
\textit{\textbf{Myopic-CE}} & - & 88.58 $\pm$ 1.32 & +9.38\\
\textit{\textbf{Myopic-CF}} & - & 90.77 $\pm$ 2.18 & +7.19\\
\bottomrule
\end{tabular}
}
\end{table}

In general, we observe that the \NeurADP~policy consistently outperforms all variations of benchmark policies for both the \DRL~and \Myopic~cases. This can be attributed to its enhanced ability to efficiently match batches of orders with available couriers, which we explore in more detail in further experiments below. Additionally, we find that the \Myopic~policies exhibit consistently superior performance compared to the \DRL-based policies. One possible explanation for this trend is that the \Myopic~policy settings prioritize maximizing the number of incoming orders batched together when matching them with a given courier. On the other hand, the \DRL~policy tends to make simpler accept/reject decisions for each individual order and matches them individually as well (i.e., \DRL-based approach does not try to identify the best batch of orders). Furthermore, we note that among the various benchmark policy variations, `\DRL\texttt{-DC}' and `\Myopic\texttt{-DC}' consistently yield the best performance. As a result, we utilize `\DRL\texttt{-DC}' as our \DRL~policy and `\Myopic\texttt{-DC}' as our \Myopic~policy for the remainder of the experiments.

\subsubsection{Impact of Number of Couriers}
We examine the impact of the number of available couriers on the number of fulfilled orders throughout the 24-hour problem horizon. The results of these experiments are presented in Table~\ref{table:couriers_table} for 10, 15, and 20 couriers. We once again utilize the \NeurADPFixed~ceiling, whose values of average orders served are shown under ``Fixed Ceiling'', and present the percentage of orders each policy fulfills based upon this ceiling in columns labeled ``\% NeurADP Filled'', ``\% Myopic Filled'', and ``\% DRL Filled''. As before, we calculate the average number of orders served by each policy and divide it by the ceiling value. This result is then multiplied by 100, and a standard deviation is provided for each policy. The final two columns illustrate the percentage increase in the average number of orders fulfilled by the \NeurADP~policy compared to that of the \Myopic~and \DRL~policies, respectively. This metric is calculated by subtracting the average number of orders fulfilled by the benchmark policies from the average number of orders fulfilled by the \NeurADP~policy, then dividing the result by the ``Fixed Ceiling'' value, multiplied by 100. These column definitions remain consistent throughout the subsequent tables presented.
\setlength{\tabcolsep}{4.5pt}
\renewcommand{\arraystretch}{1.15}
\begin{table}[!ht]
\centering
\caption{Impact of number of couriers on order fulfillment for the Brooklyn dataset (avg. number of orders fulfilled over 20-day test window is reported for the Fixed Ceiling; performance of other settings are w.r.t. Fixed Ceiling, provided as mean$\pm$stdev).}
\label{table:couriers_table}
\resizebox{0.925\textwidth}{!}{
\begin{tabular}{P{0.125\textwidth}L{0.12\textwidth}L{0.22\textwidth}C{0.18\textwidth}C{0.15\textwidth}L{0.15\textwidth}L{0.12\textwidth}}
\toprule
\textit{\textbf{Number of Couriers}} & \textit{\textbf{Fixed Ceiling}} & \textit{\textbf{\% NeurADP Filled}} & \textit{\textbf{\% Myopic Filled}} & \textit{\textbf{\% DRL Filled}} & \textit{\textbf{\% Incr. Over Myopic}} & \textit{\textbf{\% Incr. Over DRL}}\\
\midrule
\textit{\textbf{10 couriers}} & 691.85 & 97.89 $\pm$ 1.43 & 86.37 $\pm$ 1.06 & 82.14 $\pm$ 1.52 & +11.52 & +15.75\\
\textit{\textbf{15 couriers}} & 955.00 & 97.96 $\pm$ 1.97 & 91.27 $\pm$ 1.64 & 86.04 $\pm$ 1.53 & +6.69 & +11.92\\
\textit{\textbf{20 couriers}} & 1134.70 & 98.30 $\pm$ 1.98 & 95.21 $\pm$ 2.23 & 90.22 $\pm$ 1.79 & +3.09 & +8.08\\
\bottomrule
\end{tabular}
}
\end{table}

\begin{figure}[!ht]
    \centering
    \subfloat[Average return time (in minutes) of couriers making deliveries.\centering]{{\includegraphics[width=0.45\textwidth]{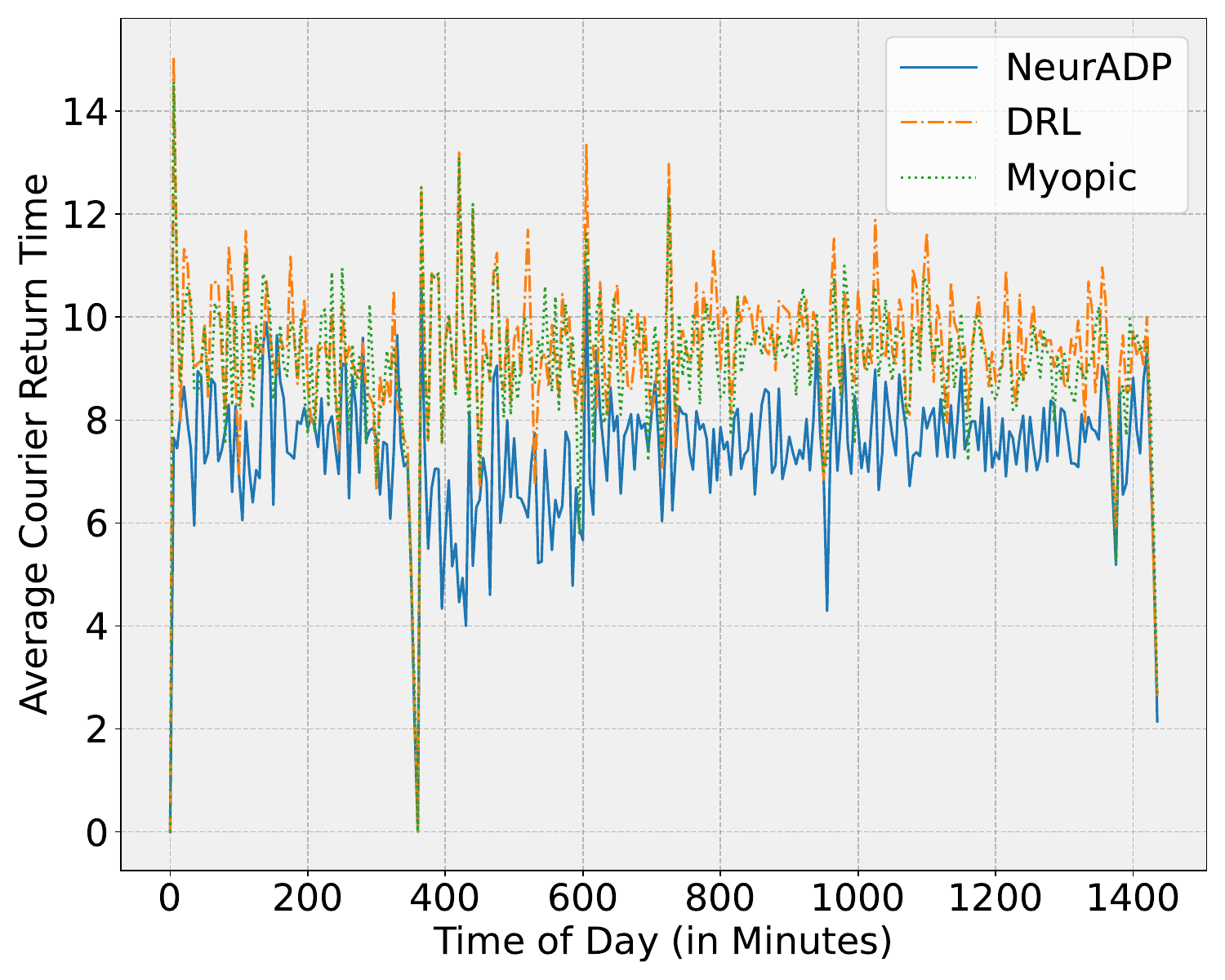} }}
    \subfloat[Average travel time (in minutes) of orders accepted. \centering]{{\includegraphics[width=0.45\textwidth]{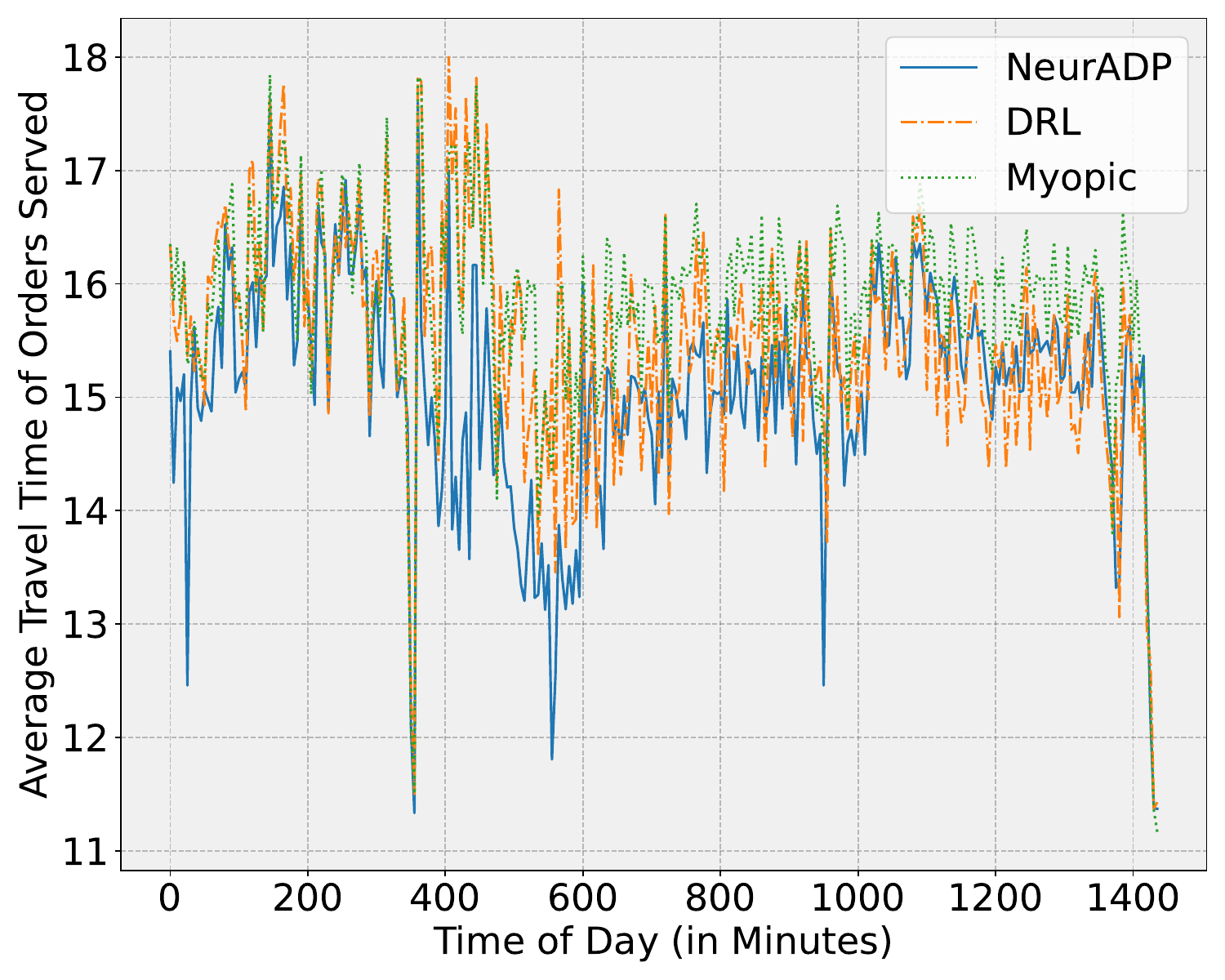} }}%
    \caption{Brooklyn dataset base-case auxiliary statistics.}
    \label{fig:couriers_stats}
\end{figure}

We observe that, once again, \NeurADP~is able to consistently outperform the benchmark policies. This superiority in performance is due to several factors. First, the \NeurADP~policy enables, on average, a greater number of agents to be available at the warehouse, compared to the benchmark policies. More specifically, in the scenario with 10 couriers, there are on average 0.37 couriers available at the warehouse at each decision epoch throughout the day, while there are only 0.06 for each of the \DRL~and \Myopic~policies. This is due to \NeurADP's ability to match batches of orders to couriers which minimize their travel time away from the warehouse, yet maximize their orders fulfilled. More specifically, unlike the other policies which look to fill the queue of each courier at each time step, the \NeurADP~policy looks to rather match couriers with orders which best allow them to return to the warehouse as soon as possible. Doing so allows them to be able to be available to more incoming orders in subsequent time steps. This can be better seen in Figure~\ref{fig:couriers_stats}, in which we see that the average return time of a courier making a delivery is consistently lower than both benchmark policies for the \NeurADP. Additionally, we see that \NeurADP~accepts orders which on average require less travel time to be delivered. More specifically, the average direct travel time from the depot to the delivery location for an accepted order is 14.41 minutes for \NeurADP, while it is 15.24 and 15.61 minutes for the \DRL~and \Myopic~policies, respectively. Furthermore, the average queue size of couriers making deliveries is 1.42 orders for \NeurADP, while it is respectively 1.84 and 2.02 for the \DRL~and \Myopic~cases. Thus, by ignoring orders at each time step, which may be out of the way for couriers, and by prioritizing return times to the depot, \NeurADP~is able to more efficiently serve orders throughout the day. This superiority in performance deteriorates, however, as more couriers are incorporated into the environment. This can be attributed to the diminishing significance of the quality of the employed policy as more couriers become available. In other words, when there are a large number of couriers available to handle order deliveries at each time step, the specific policy being used becomes less critical since most policies are capable of performing relatively well in such scenarios.

\subsubsection{Impact of Delay Time}
We next assess the impact of the permitted delay time on the number of orders served. The summary results for this experiment are provided in Table~\ref{table:delay_table}. We employ the \Direct~ceiling to establish the upper bound on the number of orders fulfilled and illustrate the percentage of fulfilled orders by each policy based on this ceiling. Similar to previous experiments, we observe that the \NeurADP~policy outperforms the two benchmark policies for all variations of delay time. Moreover, we once again notice a decline in relative performance improvements attributed to the \NeurADP~policy compared to the two benchmark policies as the problem setting becomes less restrictive (i.e., when the maximum allowed delay time is increased). 
\setlength{\tabcolsep}{4.5pt}
\renewcommand{\arraystretch}{1.15}
\begin{table}[!ht]
\centering
\caption{Impact of delay time on order fulfillment for the Brooklyn dataset (avg. number of orders fulfilled over 20-day test window is reported for the Direct Ceiling; performance of other settings are w.r.t. Direct Ceiling, provided as mean$\pm$stdev).}
\label{table:delay_table}
\resizebox{0.995\textwidth}{!}{
\begin{tabular}{P{0.125\textwidth}L{0.12\textwidth}L{0.22\textwidth}C{0.18\textwidth}C{0.15\textwidth}L{0.15\textwidth}L{0.12\textwidth}}
\toprule
\textit{\textbf{Delay Time}} & \textit{\textbf{Direct Ceiling}} & \textit{\textbf{\% NeurADP Filled}} & \textit{\textbf{\% Myopic Filled}} & \textit{\textbf{\% DRL Filled}} & \textit{\textbf{\% Incr. Over Myopic}} & \textit{\textbf{\% Incr. Over DRL}}\\
\midrule
\textit{\textbf{5 minutes}} & 1560.75 & 55.80 $\pm$ 0.95 & 49.95 $\pm$ 0.84 & 46.33 $\pm$ 0.76 & +5.85 & +9.47\\
\textit{\textbf{10 minutes}} & 1578.65 & 59.26 $\pm$ 0.62 & 55.21 $\pm$ 0.92 & 52.04 $\pm$ 0.50 & +4.04 & +7.22\\
\textit{\textbf{15 minutes}} & 1584.75 & 59.65 $\pm$ 0.50 & 56.84 $\pm$ 0.47 & 53.51 $\pm$ 0.47 & +2.81 & +6.14\\
\bottomrule
\end{tabular}
}
\end{table}

\begin{figure}[!ht]
\begin{center}
\includegraphics[width=0.60\textwidth]{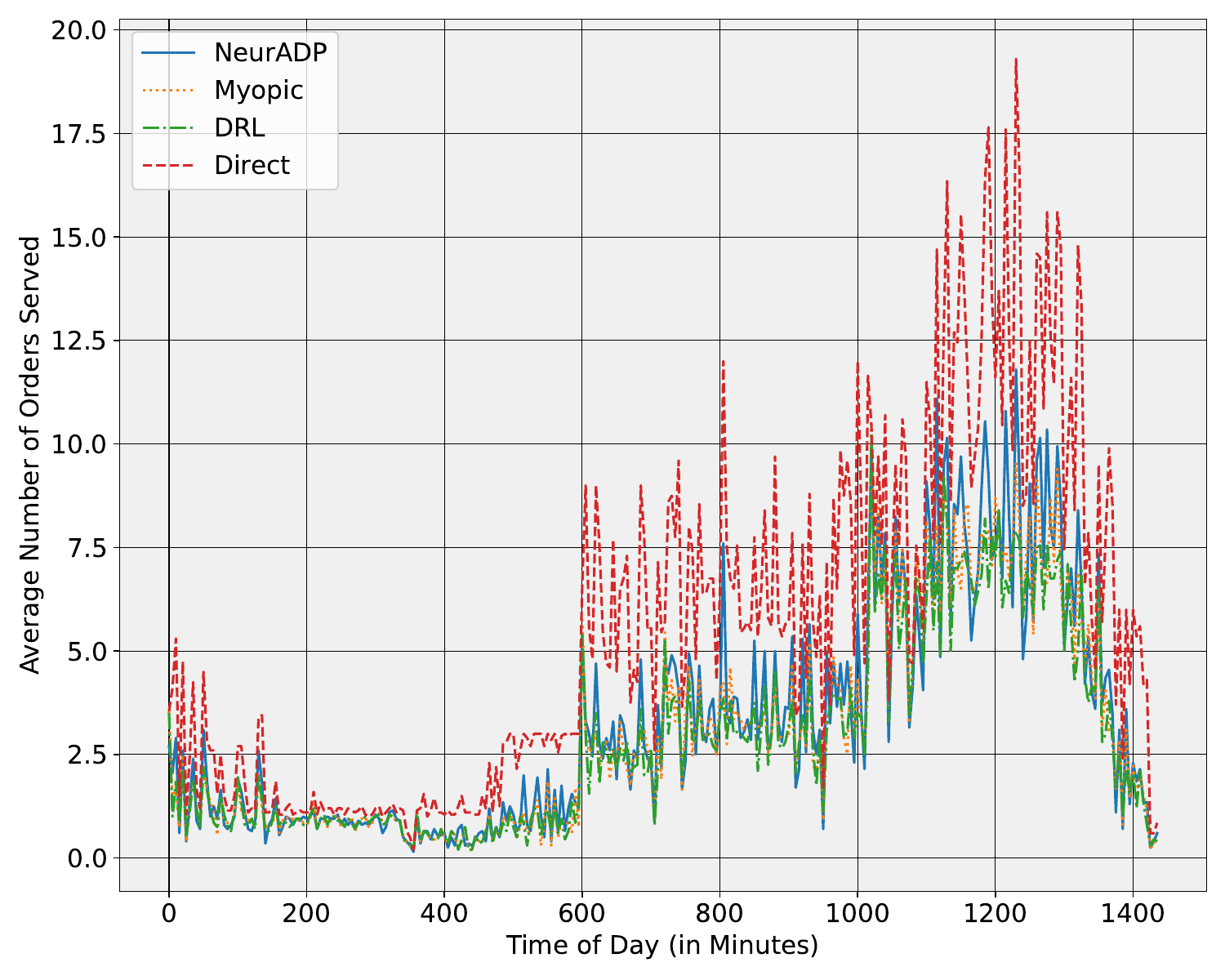}
\end{center}
\caption{Orders seen and fulfilled throughout the day for the Brooklyn dataset.}
\label{fig:orders_fulfilled}
\end{figure}

Figure~\ref{fig:orders_fulfilled} illustrates the fulfillment of orders by each policy over the problem horizon in the 10-minute delay scenario. Initially, when there are relatively few incoming orders, all policies exhibit similar performance in fulfilling the incoming orders. However, as the day progresses and a higher volume of orders arrives at each decision epoch, the \NeurADP~policy surpasses the benchmark policies, especially during peak hours. This observation highlights that when the number of orders is low and there are sufficient couriers available, all policies perform relatively well. However, as the number of orders increases and couriers become busier, making informed and intelligent decisions becomes crucial. In such situations involving ultra-fast delivery during peak periods, the \NeurADP~policy demonstrates superior performance, effectively outperforming the other policies.

\subsubsection{Impact of Courier Capacity}
Table~\ref{table:capacity_table} presents the results on the impact of the maximum capacity size of a courier and the number of orders fulfilled. These results reaffirm the superiority of the \NeurADP~policy compared to the benchmark policies across different capacity variations. Interestingly, the findings suggest that the addition of extra capacity yields diminishing returns. Specifically, increasing the capacity from 1 to 2 orders results in an average increase of 172.51 orders served among the three policies. However, the subsequent increases in capacity from 2 to 3 orders and from 3 to 4 orders correspond to smaller increases of 107.67 and 64.68 orders served, respectively. This trend is also observed in the average filled queue size of couriers making deliveries. The increase in filled queue size from a maximum capacity of 1 to 2 orders is 0.44 orders on average, while the increases from 2 to 3 orders and 3 to 4 orders are only 0.27 and 0.11 orders, respectively. It is important to note that capacity depends on various factors, such as delay time, and increasing capacity may not necessarily translate into serving more orders due to time constraints or couriers available. Thus, our findings suggest that while the \NeurADP~policy exhibits improvement over the \DRL~policy up to a capacity of 3 orders, the improvement diminishes for a capacity of 4 orders. This indicates that as the capacity increases, the problem setting becomes less restrictive, allowing for a wider range of effective policies, which benefits the \NeurADP~policy and enables it to benefit from smarter decision-making. However, as the capacity continues to increase, the improvement of \NeurADP~over the \DRL~policy declines, suggesting that well-performing policies can be more easily obtained and that the quality of the policy becomes less crucial when the problem setting becomes too non-restrictive. The improvement of \NeurADP~over the \Myopic~policy follows a similar general trend as well.
\setlength{\tabcolsep}{4.5pt}
\renewcommand{\arraystretch}{1.15}
\begin{table}[!ht]
\centering
\caption{Impact of courier capacity on order fulfillment for the Brooklyn dataset (avg. number of orders fulfilled over 20-day test window is reported for the Direct Ceiling; performance of other settings are w.r.t. Direct Ceiling, provided as mean$\pm$stdev).}
\label{table:capacity_table}
\resizebox{0.995\textwidth}{!}{
\begin{tabular}{P{0.125\textwidth}L{0.12\textwidth}L{0.22\textwidth}C{0.18\textwidth}C{0.15\textwidth}L{0.15\textwidth}L{0.12\textwidth}}
\toprule
\textit{\textbf{Courier Capacity}} & \textit{\textbf{Direct Ceiling}} & \textit{\textbf{\% NeurADP Filled}} & \textit{\textbf{\% Myopic Filled}} & \textit{\textbf{\% DRL Filled}} & \textit{\textbf{\% Incr. Over Myopic}} & \textit{\textbf{\% Incr. Over DRL}}\\
\midrule
\textit{\textbf{1 order}} & 993.10 & 62.82 $\pm$ 0.96 & 58.75 $\pm$ 1.01 & 58.48 $\pm$ 1.05 & +4.07 & +4.34\\
\textit{\textbf{2 orders}} & 1481.95 & 55.03 $\pm$ 0.69 & 51.47 $\pm$ 0.65 & 49.08 $\pm$ 0.43 & +3.56 & +5.95\\
\textit{\textbf{3 orders}} & 1578.65 & 59.26 $\pm$ 0.62 & 55.21 $\pm$ 0.92 & 52.04 $\pm$ 0.50 & +4.04 & +7.22\\
\textit{\textbf{4 orders}} & 1614.00 & 61.10 $\pm$ 0.55 & 58.54 $\pm$ 0.60 & 55.25 $\pm$ 0.45 & +2.56 & +5.85\\
\bottomrule
\end{tabular}
}
\end{table}

\subsubsection{Impact of Geographic Location}
The experimental results related to different geographic locations are presented in Table~\ref{table: locations_table}, which include Chicago, Brooklyn, Iowa, and Bangalore datasets. Each dataset has its unique distribution of order and delivery locations, as depicted in Figure~\ref{fig:data_boxplot} and Figure~\ref{fig:data_summary}. We observe that in datasets where the delivery area is more concentrated, like the Chicago dataset, all policies perform well in fulfilling order requests, and the performance advantage of the \NeurADP~policy over the benchmark policies is relatively small. However, as the delivery area expands and the delivery locations become more scattered, the improvement of the \NeurADP~policy over the benchmark policies becomes more significant. This can be explained by the following reasoning: in dense areas where delivery locations are close to the depot, making smarter decisions in matching couriers to orders or rejecting certain orders to wait for the next time step becomes less critical. This is because even with sub-optimal matching, the courier will still be able to fulfill all assigned orders and return to the warehouse on time for the next time step. However, in datasets with sparse delivery locations, such as the Bangalore dataset, where travel time from the depot to each location is longer, making sub-optimal matching decisions between couriers and incoming orders becomes more costly. A courier being occupied with sub-optimal assignments for a longer duration means they are unavailable to serve new orders, leading to delays. This becomes more apparent when comparing the average return times of couriers for each policy across the different datasets. For the Chicago dataset, the average return time for couriers making deliveries under the \NeurADP~policy is 3.63 minutes, while it is 4.63 and 4.22 minutes for the \DRL~and \Myopic~policies, respectively, showing a relatively small difference. In contrast, for the Bangalore dataset, the average return time for a courier in the \NeurADP~policy is 18.27 minutes, while it is 30.66 and 30.03 minutes for the \DRL~and \Myopic~policies, respectively, demonstrating a much larger disparity. As such, making poor matching decisions in sparse scenarios incurs higher costs, and having a smarter policy becomes significantly more important.
\setlength{\tabcolsep}{4.5pt}
\renewcommand{\arraystretch}{1.15}
\begin{table}[!ht]
\centering
\caption{Impact of geographic location on order fulfillment (avg. number of orders fulfilled over 20-day test window is reported for the Direct Ceiling for each location; performance of the policies are w.r.t. Direct Ceiling, provided as mean$\pm$stdev).}
\label{table: locations_table}
\resizebox{0.995\textwidth}{!}{
\begin{tabular}{P{0.125\textwidth}L{0.12\textwidth}L{0.22\textwidth}C{0.18\textwidth}C{0.15\textwidth}L{0.15\textwidth}L{0.12\textwidth}}
\toprule
\textit{\textbf{Geographic Location}} & \textit{\textbf{Direct Ceiling}} & \textit{\textbf{\% NeurADP Filled}} & \textit{\textbf{\% Myopic Filled}} & \textit{\textbf{\% DRL Filled}} & \textit{\textbf{\% Incr. Over Myopic}} & \textit{\textbf{\% Incr. Over DRL}}\\
\midrule
\textit{\textbf{Chicago}} & 1601.60 & 96.93 $\pm$ 0.37 & 96.61 $\pm$ 0.32 & 94.51 $\pm$ 0.35 & +0.33 & +2.43\\
\textit{\textbf{Brooklyn}} & 1578.65 & 59.26 $\pm$ 0.62 & 55.21 $\pm$ 0.92 & 52.04 $\pm$ 0.50 & +4.04 & +7.22\\
\textit{\textbf{Iowa}} & 1549.90 & 38.19 $\pm$ 0.46 & 30.63 $\pm$ 0.50 & 30.02 $\pm$ 0.61 & +7.56 & +8.17\\
\textit{\textbf{Bangalore}} & 1468.80 & 24.14 $\pm$ 0.24 & 16.88 $\pm$ 0.33 & 15.54 $\pm$ 0.52 & +7.25 & +8.60\\
\bottomrule
\end{tabular}
}
\end{table}
\subsection{Computational Performance of NeurADP}
Lastly, we provide a comprehensive overview of the auxiliary statistics related to the NeurADP algorithm drawn from our experiments. 
Notably, the algorithm takes approximately 8 hours to execute in our baseline experiment. Within a NeurADP iteration, the computational time is predominantly consumed by three tasks: data preparation for the neural network and evaluation of feasible actions (42.48\%), the generation of feasible actions while satisfying the constraints (28.83\%), and the \texttt{MatchingIP}-driven action selection process (26.39\%). Other tasks collectively account for less than 1\% of the total computational time. The considerable time spent on generating feasible actions highlights the value of experience sampling, suggesting that storing such actions for future reference could lead to significant computational savings. Our results also show that including auxiliary information in our post-decision state gives a modest performance enhancement of 0.61\%. Considering the less pronounced geographic competition among couriers in the ODP compared to the ride-pool matching problem, it can be expected that such additional information about other system agents might not yield significant benefits. Nevertheless, given the slight performance boost and minimal computational overhead, we have opted to incorporate this auxiliary information in our experiments.

\section{Conclusion}\label{conclusion}
This paper addresses the challenges and complexities of the same-day delivery problem by focusing on the order dispatching and matching aspects. Our work builds upon existing research and contributes to the literature by introducing innovative features and capabilities. It proposes the incorporation of batching and courier queues to enhance dispatching operations, providing a more realistic representation of the order dispatching process. Additionally, the scope of the problem is expanded to consider larger problem sizes, capturing the complexities of managing larger-scale dispatching operations. Furthermore, our paper introduces the application of the NeurADP approach to solving ultra-fast ODP, extending the potential applications of NeurADP beyond its original context. The effectiveness of NeurADP is demonstrated through implementation and comparison with myopic and DRL baselines, highlighting its advantages. Original datasets tailored for order dispatching operations are introduced to support the research and facilitate comprehensive evaluations. The paper also presents artificial bounds for evaluating solution quality and conducts a sensitivity analysis to investigate the performance of NeurADP under various factors. Overall, this work contributes to advancing the understanding and applicability of solution methodologies in the field of order dispatching and same-day delivery, providing valuable insights for practitioners and future research endeavors.

Future work may aim to enhance the representation of uncertainty in this problem setting by introducing loading and pickup delays, thereby capturing real-world dynamics more accurately. Additionally, exploring the incorporation of time-dependent uncertainty in order arrivals would provide insights into how service efficiency and quality impact the frequency of incoming order requests. Furthermore, multi-modal delivery within this scenario, encompassing various transportation modes like drones, autonomous vehicles, and traditional couriers can be investigated in future works as well. This approach would address the emerging trend of utilizing diverse delivery methods to optimize efficiency and address a range of delivery scenarios. Finally, fairness aspects of delivery services in this expedited delivery setting, aiming to ensure equitable access to timely deliveries across all geographical areas constitute an interesting research avenue.

\section*{Disclosure statement}
No potential conflict of interest was reported by the authors.

\section*{Data Availability Statement}
Full research data can be accessed through following link: \\ \url{https://tinyurl.com/yc32pwrp}


\singlespacing
\bibliographystyle{elsarticle-harv} 
\bibliography{ADP_odp}

\end{document}